\documentclass[a4paper, 10pt]{article}
\usepackage{amsmath}
\usepackage{amssymb,esint}
\usepackage{amscd}
\usepackage{xspace}
\usepackage{fancyhdr}
\usepackage{xcolor}
\newtheorem{theorem}{Theorem}[section]

\newtheorem{definition}[theorem]{Definition}

\newtheorem{lemma}[theorem]{Lemma}

\newtheorem{proposition}[theorem]{Proposition}
\newtheorem{remark}[theorem]{Remark}

\newenvironment{proof}[1][Proof]{\textbf{#1.} }{\hfill\rule{0.5em}{0.5em}}
{\catcode`\@=11\global\let\AddToReset=\@addtoreset
\AddToReset{equation}{section}

\AddToReset{theorem}{section}

\numberwithin{equation}{section}

\newcommand{\be}{\begin{equation}}
\newcommand{\bel}[1]{\begin{equation}\label{#1}}
\newcommand{\ee}{\end{equation}}
\newcommand{\barr}{\begin{eqnarray}}
\newcommand{\earr}{\end{eqnarray}}
\newcommand{\bars}{\begin{eqnarray*}}
\newcommand{\ears}{\end{eqnarray*}}

\newtheorem{subn}{\name}

\newcommand{\bsn}[1]{\def\name{#1}\begin{subn}}
\newcommand{\esn}{\end{subn}}
\newtheorem{sub}{\name}[section]
\newcommand{\bs}{\begin{sub}}
\newcommand{\es}{\end{sub}}
\newcommand{\bsl}[1]{\begin{sub}\label{#1}}
\newcommand{\bth}[1]{\def\name{Theorem}
\begin{sub}\label{t:#1}}
\newcommand{\blemma}[1]{\def\name{Lemma}
\begin{sub}\label{l:#1}}
\newcommand{\bcor}[1]{\def\name{Corollary}
\begin{sub}\label{c:#1}}
\newcommand{\bdef}[1]{\def\name{Definition}
\begin{sub}\label{d:#1}}
\newcommand{\bprop}[1]{\def\name{Proposition}
\begin{sub}\label{p:#1}}

\newcommand{\BA}{\begin{array}}
\newcommand{\EA}{\end{array}}
\newcommand{\BAN}{\renewcommand{\arraystretch}{1.2}
\setlength{\arraycolsep}{2pt}\begin{array}}
\newcommand{\BAV}[2]{\renewcommand{\arraystretch}{#1}
\setlength{\arraycolsep}{#2}\begin{array}}
\newcommand{\BSA}{\begin{subarray}}
\newcommand{\ESA}{\end{subarray}}
\newcommand{\BAL}{\begin{aligned}}
\newcommand{\EAL}{\end{aligned}}
\newcommand{\BALG}{\begin{alignat}}
\newcommand{\EALG}{\end{alignat}}
\newcommand{\BALGN}{\begin{alignat*}}
\newcommand{\EALGN}{\end{alignat*}}

\newcommand{\qeda}{\hspace{10mm}\hfill $\square$}

\newcommand{\abs}[1]{\left |#1\right |}

\def\angb<#1>{\langle #1 \rangle}

\newcommand{\opname}[1]{\mbox{\rm #1}\,}

\newcommand{\dist}{\opname{dist}}
\newcommand{\myfrac}[2]{{\displaystyle \frac{#1}{#2} }}
\newcommand{\myint}[2]{{\displaystyle \int_{#1}^{#2}}}

\newcommand{\prt}{\partial}

\newcommand{\nind}{\noindent}

\def\ga{\alpha}     \def\gb{\beta}       \def\gg{\gamma}
       \def\gd{\delta}      \def\ge{\epsilon}
\def\gth{\theta}                         
           
      \def\gk{\kappa}      \def\gl{\lambda}
\def\gm{\mu}        \def\gn{\nu}         
    \def\gr{\rho}        
       
      \def\gw{\omega}
                
     \def\Gd{\Delta}

\def\Gw{\Omega}              

   \def\BBN {\mathbb N}    
   \def\BBR {\mathbb R}

\begin{document}
\title{Quasilinear and Hessian Lane-Emden type\\
	systems with measure data
}
\author{
{\bf Marie-Fran\c{c}oise Bidaut-V\'eron\thanks{E-mail address: veronmf@univ-tours.fr, Laboratoire de Math\'ematiques et Physique Th\'eorique, Universit\'e Fran\c{c}ois Rabelais,  Tours,  France}}\\[0.5mm]
{\bf Quoc-Hung Nguyen\thanks{ E-mail address: quochung.nguyen@sns.it,  Scuola Normale Superiore, Centro Ennio de Giorgi, Piazza dei Cavalieri 3, I-56100
		Pisa, Italy.}}\\[0.5mm]
{\bf Laurent V\'eron\thanks{ E-mail address: veronl@lmpt.univ-tours.fr, Laboratoire de Math\'ematiques et Physique Th\'eorique, Universit\'e Fran\c{c}ois Rabelais,  Tours,  France}}\\[2mm]}
\date{}  
\maketitle
\tableofcontents
\begin{abstract}  We study nonlinear systems of the form $-\Delta_pu=v^{q_1}+\mu,\;
-\Delta_pv=u^{q_2}+\eta$ and  $F_k[-u]=v^{s_1}+\mu,\;
F_k[-v]=u^{s_2}+\eta$ in a bounded domain $\Omega$ or in $\mathbb{R}^N$ where $\mu$ and $\eta$ are nonnegative Radon measures, $\Delta_p$ and $F_k$ are respectively the $p$-Laplacian and the $k$-Hessian operators and $q_1$, $q_2$, 
$s_1$ and $s_2$ positive numbers. We give necessary and sufficient conditions for existence expressed in terms of Riesz or Bessel capacities.
\vspace{2mm}
\noindent {\it \footnotesize 2010 Mathematics Subject Classification}. {\scriptsize
	35J70, 35J60, 45G15, 31C15}.

\noindent {\it \footnotesize Key words:} {\scriptsize $p$-Laplacian, $k$-Hessian, Bessel and Riesz capacities, measures, maximal functions.}

\end{abstract}
\section{Introduction and Main results}

Let $\Omega\subset\mathbb{R}^N$ be either a bounded domain or the whole $\mathbb{R}^N$, $p>1$ and $k\in\{1,2,...,N\}$. We denote by 
$$\Delta_pu:=div\left(|{\nabla u}|^{p-2} \nabla u\right)$$ 
the p-Laplace operator and by 
$$ F_k[u]=\sum_{1\leq j_1<j_2<...<j_k\leq N}\lambda_{ j_1}\lambda_{ j_2}...\lambda_{ j_k}$$
the k-Hessian operator where $\lambda_1,...,\lambda_N$ are the eigenvalues of the Hessian matrix $D^2u$. In the work \cite {22PhVe}, Phuc and Verbitsky  obtained necessary and sufficient conditions for existence of nonnegative solutions to the following equations 
\begin{align}\label{2hvEQ1*}
\begin{array}{ll}
-\Delta_pu=u^q+\mu\qquad&\text{in }\Omega\\
\phantom{ -\Delta_p}
u=0\qquad&\text{on }\partial\Omega,
\end{array}
\end{align}
and 
\begin{align}
\label{2hvEQ2*}
\begin{array}{ll}
F_k[-u]=u^q+\mu\qquad&\text{in }\Omega\\
\phantom{  F_k[-]}
u= 0\qquad&\text{on }\partial\Omega.
\end{array}
\end{align}
Their conditions  involve the continuity of the measures with respect to Bessel or Riesz capacities and Wolff potentials estimates. 
For example, if $\Omega$ is bounded and $\mu$ has compact support in $\Omega$, they proved that it is equivalent to solve  
(\ref{2hvEQ1*}), or to have
\begin{equation}
\label{2hvEQ3}\begin{array}{ll}
\mu(E)\leq c_1\text{Cap}_{\mathbf{G}_p,\frac{q}{q+1-p}}(E)\qquad\text{for all compact set }E\subset\Omega,
\end{array}
\end{equation}
for some constant  $c_1>0$ where $\text{Cap}_{\mathbf{G}_p,\frac{q}{q+1-p}}$ is a Bessel capacity, or to have
\begin{equation}
\label{2hvEQ4}\begin{array}{ll}
\myint{B}{}\left({\bf W}^{R}_{1,p}[\mu_B](x)\right)^{q}dx\leq c_2\mu (B)\qquad\text{for all ball }B \text{ s.t. }
B\cap\text{supp}\mu\not=\emptyset,
\end{array}
\end{equation}
for some constant  $c_2>0$, where $R=2\text{ diam}(\Omega)$ and ${\bf W}^{R}_{1,p}[\mu_B]$ denotes the $R$-truncated Wolff potential of the measure $\mu_B=\chi_{_B}\gm$.  Concerning the k-Hessian operator in a bounded $(k-1)$-convex domain $\Omega$,  they proved that if $\mu$ has compact support,  the  problem \eqref{2hvEQ2*} with $q>k$ admits a nonnegative solution if and only if 
\begin{equation}
\label{2hvEQ6}\mu(E)\leq c_3\text{Cap}_{\mathbf{G}_{2k},\frac{q}{q-k}}(E)\qquad\text{for all compact set }E\subset\Omega,
\end{equation}
for some $c_3$. In turn this condition is equivalent to 
\begin{equation}
\label{2hvEQ7}\int_B\left[{\bf W}^{R}_{\frac{2k}{k+1},k+1}[\mu_B(x)]\right]^{q}dx\leq c_4\mu (B)\qquad\text{for all ball }B \text{ s.t. }
B\cap\text{supp}\,\mu\not=\emptyset,
\end{equation}
for some $c_4>0$. The results concerning the  linear case $p=2$ and $k=1$, can be found in \cite{22AdPi,22BaPi,22V}.\medskip\\
The natural counterpart of equation \eqref{2hvEQ1*} and \eqref{2hvEQ2*} for  systems:
\begin{align}\label{2hvEQ1}
\begin{array}{ll}
-\Delta_pu=v^{q_1}+\mu\qquad&\text{in }\Omega\\
-\Delta_pv=u^{q_2}+\eta\qquad&\text{in }\Omega\\
\phantom{ -\Delta_p}
u=v=0\qquad&\text{on }\partial\Omega,
\end{array}
\end{align}
and 
\begin{align}
\label{2hvEQ2}
\begin{array}{ll}
F_k[-u]=v^{s_1}+\mu\qquad&\text{in }\Omega\\
F_k[-v]=u^{s_2}+\eta\qquad&\text{in }\Omega\\\phantom{  F_k[-]}
u= v=0\qquad&\text{on }\partial\Omega,
\end{array}
\end{align}
where $q_1,q_2>p-1,s_1,s_2>k$ and $\mu,\eta$ are Radon measures. If $\Omega=\mathbb{R}^N$, we consider the same equations, except that the boundary conditions are replaced by $\inf_{\mathbb{R}^N}u=\inf_{\mathbb{R}^N}v=0$ and  our statements involve the Riesz potentials and their associated capacities $\operatorname{Cap}_{I_{\ga,\gb}}$.
Our main results are the following.  \medskip

\nind{\bf Theorem A} {\it
	 Let $1<p<N$, $q_1,q_2>0$ and $q_2q_1>(p-1)^2$.  Let $\mu,\eta$ be  nonnegative Radon measures in $\mathbb{R}^N$. If 
	 the following system
	\begin{equation}\label{2hvMT1ab}
	\left. \begin{array}{ll}
	- \Delta _p u = v^{q_1}+\mu \qquad&\text{in }\;\mathbb{R}^N  \\ 
	- \Delta _p v = u^{q_2}+\eta \qquad&\text{in }\;\mathbb{R}^N,  \\ 
	\end{array} \right.
	\end{equation} 
	admits  a nonnegative p-superharmonic solution $(u,v)$ then there exists a positive constant $c_5$ depending on $N,p,q_1,q_2$  such that 
	\begin{equation}\label{2hvMT1cb}
		\eta(E)+\int_E \left(\mathbf{W}_{1,p}[\mu](x)\right)^{q_2}dx\leq c_5 \operatorname{Cap}_{\mathbf{I}_{\frac{p(q_1+p-1)}{q_1},\frac{q_1q_2}{q_1q_2-(p-1)^2}}}(E)~~\textrm{ for all Borel sets } E.
	\end{equation}
	 Conversely, if $\gm$  and $\eta$ are bounded, there exists $c_6>0$ depending on $N,p,q_1,q_2$ such that if  $0< q_1<\frac{N(p-1)}{N-p}$ and  \eqref{2hvMT1cb} holds with $c_5$ replaced by  $c_6$, then  \eqref{2hvMT1ab} admits a nonnegative p-superharmonic solution $(u,v)$ satisfying  
	\begin{align}\label{thm1-es-uv}
	v\leq c_8 \mathbf{W}_{1,p}[\omega],~~~u\leq c_9 \mathbf{W}_{1,p}[\left(\mathbf{W}_{1,p}[\omega]\right)^{q_1}]+c_7 \mathbf{W}_{1,p}[\mu]
	\end{align} 
	in $\BBR^N$ for some $c_7,c_8,c_9>0$ where $d\omega=\left(\mathbf{W}_{1,p}[\mu]\right)^{q_2}dx+d\eta$  .
}\medskip

We notice that the left-hand side in (\ref{2hvMT1cb}) is not symmetric in $\eta$ and $\gm$ and the capacity in the right-hade side is not symmetric in $q_1$ and $q_2$. Hence the following symmetrized inequality holds
	\begin{equation}\label{2hvMT1cb+1}
		\gm(E)+\int_E \left(\mathbf{W}_{1,p}[\eta](x)\right)^{q_1}dx\leq c'_5 \operatorname{Cap}_{\mathbf{I}_{\frac{p(q_2+p-1)}{q_2},\frac{q_1q_2}{q_1q_2-(p-1)^2}}}\!\!\!\!\!\!(E)~~\textrm{ for all Borel sets } E.
	\end{equation}

It is known that 
\begin{align*}
\operatorname{Cap}_{\mathbf{I}_{\alpha,\beta}}(K)=0 ~~~\forall K \text{ compact},
\end{align*}
if $\alpha\beta\geq N$, the first part of above implies the following Liouville theorem, obtained by another method in \cite[Th 5.3\,-(i)]{22Bi4}. 
\medskip

\nind{\bf Corollary B} {\it Assume that 
	\begin{align*}
	\frac{p(q_1q_2+(p-1)\max\{q_1,q_2\})}{q_1q_2-(p-1)^2}\geq N.
	\end{align*} Any nonnegative p-superharmonic solution $(u,v)$ of inequalities
	\begin{equation}\label{2hvMT1ab*}
	\left. \begin{array}{ll}
	- \Delta _p u \geq  v^{q_1}\qquad&\text{in }\;\mathbb{R}^N  \\ 
	- \Delta _p v \geq u^{q_2} \qquad&\text{in }\;\mathbb{R}^N,  \\ 
	\end{array} \right.
	\end{equation}
	is trivial, i.e. $u=v=0$. 
}\medskip

Classical Liouville results for one equation or inequality, are proved in \cite{22Bi1}, \cite{22Bi2}, \cite{22BiDe}, \cite{22SeZo}. \medskip

When $\Omega$ is bounded domain, we have a similar result in which we denote by $d$ the distance function to the boundary $x\mapsto d(x)=\dist (x,\prt\Gw)$.\medskip

\nind{\bf Theorem C} {\it
	Let $1<p<N$, $q_1,q_2>0$ and $q_2q_1>(p-1)^2$. Let $\Omega\subset \mathbb{R}^N$ be a bounded domain and $\mu,\eta$  nonnegative Radon measures in $\Omega$.  If the following problem
	\begin{equation}\label{2hvMT1a}
	\left. \begin{array}{ll}
	- \Delta _p u = v^{q_1}+\mu \qquad&\text{in }\;\Omega  \\ 
	- \Delta _p v = u^{q_2}+\eta \qquad&\text{in }\;\Omega  \\ 
	\phantom{    - \Delta _p}
	u =v= 0 \qquad&\text{on }\;\partial\Omega,
	\end{array} \right.
	\end{equation}
	admits a nonnegative renormalized solution $(u,v)$, then then for any compact set $K\subset \Omega$, there exists a positive constant $c_{_{10}}$ depending on $N,p,q_1,q_2$ and $\text{dist}(K,\partial\Omega)$ such that 
	\begin{equation}\label{2hvMT1c}
	\eta(E)+\int_E \left(\mathbf{W}^{\frac{d(x)}{4}}_{1,p}[\mu](x)\right)^{q_2}dx\leq c_{_{10}} \operatorname{Cap}_{\mathbf{G}_{\frac{p(q_1+p-1)}{q_1},\frac{q_1q_2}{q_1q_2-(p-1)^2}}}(E)~~\textrm{ for all Borel sets } E\subset K.
	\end{equation}
	Conversely, let $\gm$  and $\eta$ be bounded with the property that there exists $c_{_{11}}>0$ depending on $N,p,q_1,q_2$ and $R=2diam\,(\Omega)$ such that if $0<q_1<\frac{N(p-1)}{N-p}$ and  \begin{align}
	 \eta (K)+\int_{K}\left(\mathbf{W}^{2R}_{1,p}[\mu]\right)^{q_2}dx\leq c_{_{11}} \operatorname{Cap}_{\mathbf{G}_{\frac{p(q_1+p-1)}{q_1},\frac{q_1q_2}{q_1q_2-(p-1)^2}}}(K),
	 \end{align}
	 for all compact set $K\subset\Gw$,  then  \eqref{2hvMT1a} admits a nonnegative renormalized solution $(u,v)$ satisfying  
	 \begin{align}\label{thm2-es-uv}
	 v\leq c_{_{13}} \mathbf{W}_{1,p}^{R}[\omega],~~~u\leq c_{_{14}} \mathbf{W}_{1,p}^{R}[\left(\mathbf{W}^{R}_{1,p}[\omega]\right)^{q_1}]+c_{_{12}} \mathbf{W}_{1,p}^{R}[\mu]
	 \end{align} 
	in $\Gw$, where  $d\omega=\left(\mathbf{W}_{1,p}^{R}[\mu]\right)^{q_2}dx+d\eta$.
}\medskip

It is known that 
\begin{align*}
\operatorname{Cap}_{\mathbf{G}_{\alpha,\beta}}(\{x_0\})>0 
\end{align*}
if and only if $\alpha\beta> N$. Thus, as an application in a partially subcritical case we have,
\medskip

\nind{\bf Corollary D} {\it Let the assumptions on $p$, $q_1$, $q_2$, $\Gw$ and $R$ of Theorem C be satisfied, $x_0\in\Gw$, $a>0$ and $\gm$ be a nonnegative Radon measures in $\Gw$. If the following problem 
\begin{equation}\label{2hvMT1a1}
	\left. \begin{array}{ll}
	- \Delta _p u = v^{q_1}+\gm \qquad&\text{in }\;\Omega  \\ 
	- \Delta _p v = u^{q_2}+a\gd_{x_0} \qquad&\text{in }\;\Omega  \\ 
	\phantom{    - \Delta _p}
	u =v= 0 \qquad&\text{on }\;\partial\Omega,
	\end{array} \right.
	\end{equation}
admits a nonnegative renormalized solution $(u,v)$, then there exist positive constants $c_{_{15}}=c_{_{15}}(N,p,q_1,q_2,d(x_0))$ and, for any compact subset $K$ of $\Gw$, $c_{_{16}}=c_{_{16}}(N,p,q_1,q_2,\dist(K,\prt\Gw) \textcolor{red}{)}$, such that
\begin{equation}
\label{sub1}
\BA {lll}
(i)\qquad\qquad\qquad\qquad &N<\myfrac{pq_2(q_1+p-1)}{q_1q_2-(p-1)^2},\qquad\qquad\qquad\qquad\qquad\qquad\\[3mm]

(ii)\qquad &a\leq c_{_{15}},\\[2mm]

 (iii) \qquad &\myint{K}{}\left(\mathbf{W}^{2R}_{1,p}[\mu]\right)^{q_2}dx\leq c_{_{16}}.
 \EA
	\end{equation}
Conversely, assuming that $\gm$ is bounded,  there exist positive constants $c_{_{17}}=c_{_{17}}(N,p,q_1,q_2,d(x_0))$, 
 $c_{_{18}}=c_{_{18}}(N,p,q_1,q_2)$ such that if $0<q_1<\frac{N(p-1)}{N-p}$ and $(\ref{sub1})$ holds with $c_{_{15}}$ and 
 $c_{_{16}}$ replaced respectively by $c_{_{17}}$ and $c_{_{18}}$, then  there exists a nonnegative renormalized solution $(u,v)$ of $(\ref{2hvMT1a1})$ satisfying 
  \begin{align}\label{thm2-es-uv2}
	 v\leq c_{_{21}} W_{1,p}^{R}[\omega],~~~u\leq c_{_{22}} \mathbf{W}_{1,p}^{R}[\left(W^{R}_{1,p}[\omega]\right)^{q_1}]+c_{_{20}} \mathbf{W}_{1,p}^{R}[\mu]
	 \end{align} 
	in $\Gw$, where  
	$$W_{1,p}^{R}[\omega]=\mathbf{W}_{1,p}^{R}\left[\left(\mathbf{W}_{1,p}^{R}[\mu]\right)^{q_2}\right]+
	a^{\frac{1}{p-1}}\left(\abs{x-x_0}^{-\frac{N-p}{p-1}}-R^{-\frac{N-p}{p-1}}\right)_+.$$
	
}\medskip


Concerning the $k$-Hessian operator we recall some notions introduced by Trudinger and Wang \cite{22TW1,22TW2,22TW3}, and we follow their notations. For $k=1,...,N$ and $u\in C^2(\Omega)$ the k-Hessian operator $F_k$ is defined by 
$$F_k[u]=S_k(\lambda(D^2u)),$$
where $\lambda(D^2u)=\lambda=(\lambda_1,\lambda_2,...,\lambda_N)$ denotes the eigenvalues of the Hessian matrix of second partial derivatives $D^2u$ and $S_k$ is the k-th elementary symmetric polynomial that is \[{S_k}(\lambda ) = \sum\limits_{1 \le {i_1} < ... < {i_k} \le N} {{\lambda _{{i_1}}}...{\lambda _{{i_k}}}}. \]
Since ${D^2}u$ is symmetric, it is clear that \[{F_k}[u] = {\left[ {{D^2}u} \right]_k},\]
where we denote by $[A]_k$  the sum of the k-th principal minors of a matrix $A=(a_{ij})$.
In order that there exists a smooth k-admissible function which vanishes on $\partial \Omega$, the boundary $\partial \Omega$ must satisfy a uniformly (k-1)-convex condition, that is 
\begin{equation*}
S_{k-1}(\kappa ) \geq c_0 >0~on ~~ \partial\Omega,
\end{equation*}
for some positive constant $c_0$, where $\kappa= (\kappa_1,\kappa_2,...,\kappa_{n-1})$ denote the principal curvatures of $\partial \Omega$ with respect to its inner normal. We also denote by $\Phi^k(\Omega)$ the class of upper-semicontinuous functions $\Omega\to[Õ-\infty,\infty)$ which are $k$-convex, or subharmonic in the Perron sense (see Definition \ref{2hvk-conv}).\
In this paper we prove the following theorem (in which expression $\mathbb E[q]$ is the largest integer less or equal to $q$)
\medskip

\nind{\bf Theorem E} {\it Let $2k<N,s_1,s_2>0$, $s_1s_2>k^2$. Let $\Omega$ be a bounded uniformly (k-1)-convex domain in $\mathbb{R}^N$ with diameter $R$. Let  $\mu=\mu_1+f$ and $\eta=\eta_1+g$ be nonnegative Radon measures where $\mu_1,\eta_1$ has compact support in $\Omega$ and $f,g\in L^l(\Omega)$ for some $l>\frac{N}{2k}$.  If the following problem 
	\begin{equation}\label{2hvMT3a}
	\left. \begin{array}{ll}
	F_k[-u] = v^{s_1}+\mu \qquad&\text {in }\;\Omega  \\ 
	F_k[-v] = u^{s_2}+\eta \qquad&\text {in }\;\Omega  \\ 
	\phantom{F_k}
	u = v=0&\text {on }\;\partial\Omega,  \\ 
	\end{array} \right.
	\end{equation}
	admits a nonnegative  solutions $(u,v)$, continuous near $\partial \Omega$, with $-u$ and $-v$ elements of $\Phi^k(\Omega)$, then for any compact set $K\subset \Omega$, there exists a positive constant $c_{_{23}}$ depending on $N,k, s_1,s_2$ and $dist(K,\partial\Omega)$ such that there holds 
	\begin{equation}\label{2hvMT3c}
	\eta(E)+\int_E \left(\mathbf{W}^{\frac{d(x)}{4}}_{\frac{2k}{k+1},k+1}[\mu](x)\right)^{s_2}dx\leq c_{_{23}} \operatorname{Cap}_{\mathbf{G}_{\frac{2k(s_1+k)}{s_1},\frac{s_1s_2}{s_1s_2-k^2}}}(E)\qquad\forall E\subset K, E\text { Borel}.
	\end{equation}
	Conversely,, if $\gm$  and $\eta$ are bounded, there exist a positive constant $c_{_{24}}$ depending on $N,k,s_1,s_2$ and $diam\,(\Omega)$ such that, if $k\leq s_1<\frac{Nk}{N-2k}$
	and
	\begin{align}
	\eta (K)+\int_{K}\left(\mathbf{W}^{2R}_{\frac{2k}{k+1},k+1}[\mu]\right)^{s_2}dx\leq c_{_{24}} \operatorname{Cap}_{\mathbf{G}_{\frac{2k(s_1+k)}{s_1},\frac{s_1s_2}{s_1s_2-k^2}}}(K)
	\end{align}
	for all Borel set $K\subset\Gw$, then \eqref{2hvMT3a} 	admits a nonnegative  solution $(u,v)$, continuous near $\partial \Omega$, with $-u,-v\in \Phi^k(\Omega)$ satisfying 
	\begin{align}\label{thm3-eswolff}
	v\leq c_{_{28}} {\bf W}^{R}_{\frac{2k}{k+1},k+1}[\omega],~~~u\leq c_{_{29}} {\bf W}^{R}_{\frac{2k}{k+1},k+1}[\left({\bf W}^{R}_{\frac{2k}{k+1},k+1}[\omega]\right)^{s_1}]+c_{_{27}} {\bf W}^{R}_{\frac{2k}{k+1},k+1}[\mu]
	\end{align}
	in $\Omega$ for some constants $c_j$ ($j=27, 28, 29$) depending on $N,k,s_1,s_2$,  and $diam\,(\Omega)$.
}\medskip

If $\Gw$ is replaced by the whole space we prove,

\medskip

\nind{\bf Theorem F} {\it Let $2k<N,s_1,s_2>0$, $s_1s_2>k^2$.  Let  $\mu,\eta$ be a nonnegative Radon measures in $\BBR^N$. If the following  problem 
	\begin{equation}\label{2hvMT3ab}
	\left. \begin{array}{ll}
	F_k[-u] = v^{s_1}+\mu \qquad&\text {in }\;\mathbb{R}^N  \\ 
	F_k[-v] = u^{s_2}+\eta \qquad&\text {in }\;\mathbb{R}^N,  \\ 
	\end{array} \right.
	\end{equation}
	admits a nonnegative  solutions $(u,v)$ with $-u$ and $-v$ belonging to $\Phi^k(\mathbb{R}^N)$, then 
	there exists a positive constant $c_{_{30}}$ depending on $N,k, s_1,s_2$ such that there holds 
	\begin{equation}\label{2hvMT3cb}
	\eta(E)+\int_E \left(\mathbf{W}_{\frac{2k}{k+1},k+1}[\mu](x)\right)^{s_2}dx\leq c_{_{30}} \operatorname{Cap}_{\mathbf{G}_{\frac{2k(s_1+k)}{s_1},\frac{s_1s_2}{s_1s_2-k^2}}}(E)\qquad\forall  E\text { Borel}.
	\end{equation}
	 Conversely,, if $\gm$  and $\eta$ are bounded,  there exists  positive constant $c_{_{31}}$ depending on $N,k,s_1,s_2$ such that, if $0< s_1<\frac{Nk}{N-2k}$ and \eqref{2hvMT3cb} holds with $c_{_{31}}$ instead of $c_{_{30}}$, then  \eqref{2hvMT3ab} admits a nonnegative solution $(u,v)$ with $-u$ and $-v$ in $\Phi^k(\mathbb{R}^N)$ satisfying 
	 	\begin{align}\label{thm4-eswolff}
	 v\leq c_{_{33}} {\bf W}_{\frac{2k}{k+1},k+1}[\omega],~~~u\leq c_{_{34}} {\bf W}_{\frac{2k}{k+1},k+1}[\left({\bf W}_{\frac{2k}{k+1},k+1}[\omega]\right)^{s_1}]+c_{_{32}} {\bf W}_{\frac{2k}{k+1},k+1}[\mu]
	 \end{align}
	 in $\mathbb{R}^N$, where the $c_{j}$ ($j=32, 33, 34$) depend on $N,k,s_1,s_2$.
}\medskip

As in the p-Laplace case, we have a Liouville property for Hessian systems.
\medskip

\nind{\bf Corollary G} {\it Assume that 
	\begin{align}
	\frac{2k(s_2s_1+k\max\{s_1,s_2\})}{s_1s_2-k^2}\geq N.
	\end{align}
Any nonnegative solution (u,v) of inequalities  
\begin{equation}\label{2hvMT3ab*}
\left. \begin{array}{ll}
F_k[-u] \geq  v^{s_1} \qquad&\text {in }\;\mathbb{R}^N  \\ 
F_k[-v] \geq  u^{s_2} \qquad&\text {in }\;\mathbb{R}^N,  \\ 
\end{array} \right.
\end{equation}
with $-u$ and $-v$ in $\Phi^k(\mathbb{R}^N)$ is trivial.
}\medskip


\section{Estimates on potentials} 
Throughout this article $c_{_j}$, j=1,2,..., denote structural positive constants and $c_{_N}$ is the volume of the unit ball in $\BBR^N$. The following inequality will be used several times in the sequel. 
\begin{lemma}\label{hardy} Let $\kappa,\gamma,\theta\in\BBR$,  such that $\kappa,\gamma>0$. Let  $h:(0,\infty)\to (0,\infty)$ be nondecreasing. Then, 
\begin{align}\label{es0}
\int_{0}^{R}t^{\kappa}\left(\int_{t}^{R}h(r)r^{\theta}\frac{dr}{r}\right)^\gamma \frac{dt}{t}\leq c_{_{35}}\int_{0}^{2R} t^{\kappa+\theta \gamma}h^\gamma (t)\frac{dt}{t}\qquad\forall R\in (0,\infty],
\end{align}
for some $c_{_{35}}>0$ depending on $\kappa$, $\gg$, $\gth$.
\end{lemma}
\proof {\it Case 1: $\gamma\leq 1$}. 
Since there holds
	\begin{align*}
	\left(\sum_{j=0}^{\infty}a_j\right)^\gamma\leq \sum_{j=0}^{\infty}a_j^\gamma \qquad\forall a_j\geq 0, 
	\end{align*}
	we deduce 
	\begin{align*}
	\left(\int_{t}^{R}h(r)r^{\theta}\frac{dr}{r}\right)^\gamma&\leq c_{_{\gg,\gth}} \left(\sum_{j=0}^{j_0}h(2^{\frac{j+1}{4}}t) (2^{\frac{j}{4}}t)^\theta \right)^\gamma\\
	&\leq c_{_{\gg,\gth}} \sum_{j=0}^{j_0}\left(h^\gamma(2^{\frac{j+1}{4}}t)\right)(2^{\frac{j}{4}}t)^{\theta\gamma} \\& \leq  c_{_{\gg,\gth}}\int_{t}^{2R}h^\gamma(r) r^{\theta\gamma}\frac{dr}{r},
	\end{align*}
	where $c_{_{\gg,\gth}}=2^{\frac{\gg}{4}}\max\{1,2^{-\frac{\gg\gth}{4}}\}$ and $2^{\frac{j_0}{4}}t<R\leq 2^{\frac{j_0+1}{4}}t$ if $R<\infty$ and $j_0=\infty$ if $R=\infty$.
	By Fubini's theorem,
	\begin{align*}
	\int_{0}^{R}t^{\kappa}\left(\int_{t}^{R}h(r)r^{\theta}\frac{dr}{r}\right)^\gamma \frac{dt}{t}&\leq c_{_{\gg,\gth}} \int_{0}^{R}t^{\kappa}\int_{t}^{2R}h^\gamma(r) r^{\theta\gamma}\frac{dr}{r} \frac{dt}{t}\\& \leq \myfrac{c_{_{\gg,\gth}}}{\gk} \int_{0}^{2R} t^{\kappa+\theta \gamma}h^\gamma(t) \frac{dt}{t},
	\end{align*}
	 which is \eqref{es0}.\medskip
	 
	 \nind{\it Case 2: $\gamma> 1$}. Since
	 $$\left(\int_{t}^{R}h(r)r^{\theta}\frac{dr}{r}\right)^\gamma 
	 \leq \left(\int_{t}^{R}r^{-\frac{\gg}{\gg-1}}\frac{dr}{r}\right)^{\gamma-1} \int_{t}^{R}h^\gg(r)r^{\gg(1+\theta)}\frac{dr}{r},
	 $$
we obtain
	\begin{align*}
		\int_{0}^{R}t^{\kappa}\left(\int_{t}^{R}h(r)r^{\theta}\frac{dr}{r}\right)^\gamma \frac{dt}{t} \leq  c_{_{\gg,\gk}}\int_{0}^{2R} t^{\kappa+\theta \gamma}h^\gamma(t) \frac{dt}{t},
	\end{align*}
by Fubini's theorem, which completes the proof.\qeda\medskip


We recall that if $\ga>0$, $1<\gb<\frac{N}{\ga}$ and $\gm$ belongs to the set of positive Radon measures in $\BBR^N$ that we denote $\mathfrak M^+(\BBR^N)$,   the Wolff potential of $\gm$ is defined by
\bel{wolff1}
{\bf W}_{\ga,\gb}[\gm](x)=\myint{0}{\infty}\left(\myfrac{\gm(B_r(x))}{r^{N-\ga\gb}}\right)^{\frac{1}{p-1}}\myfrac{dr}{r},
\ee
and if $R>0$, the $R$-truncated Wolff potential of $\gm$ is
\bel{wolff2}
{\bf W}^R_{\ga,\gb}[\gm](x)=\myint{0}{R}\left(\myfrac{\gm(B_r(x))}{r^{N-\ga\gb}}\right)^{\frac{1}{p-1}}\myfrac{dr}{r}.
\ee

If $\gm$ is a Radon measure on a Borel set $G$, it's Wolff potential (or truncated Wolff potential) is the potential of its extension by $0$ in $G^c$.
We start with the following composition estimate on Wolff potentials. 
 \begin{lemma} \label{241220142} Let $1<\beta<N/\alpha$.  Then for any $q>0$ and $\gm\in \mathfrak M^+(\BBR^N)$ we have
 	\begin{align}\label{wolffes1}
	\mathbf{W}_{\frac{\alpha\beta(q+\beta-1)}{q+(\beta-1)^2},\frac{(\beta-1)^2}{q}+1}[\mu]\leq c_{_{36}} \mathbf{W}_{\alpha,\beta}\left[\left(\mathbf{W}_{\alpha,\beta}[\mu]\right)^q\right],
 	\end{align}
in $\mathbb{R}^N$ for some $c_{_{36}}>0$ depending on $\ga,\gb,N,q$. Moreover, if $0< q <\frac{N(\beta-1)}{N-\alpha\beta}$, there holds
 	\begin{align}\label{wolffes2}
 \mathbf{W}_{\alpha,\beta}\left[\left(\mathbf{W}_{\alpha,\beta}[\mu]\right)^q\right](x)\leq c_{_{37}}	\mathbf{W}_{\frac{\alpha\beta(q+\beta-1)}{q+(\beta-1)^2},\frac{(\beta-1)^2}{q}+1}[\mu],
 	\end{align}
 	in $\mathbb{R}^N$, where $c_{_{37}}>0$ depends on $\ga,\gb,N,q$.
 \end{lemma}
 \proof For any $x\in \mathbb{R}^N$, using the fact if $y\in B_t(x)$ then $B_t(x)\subset B_{2t}(y)$, we have
 	\begin{align*}
 	\mathbf{W}_{\alpha,\beta}\left[\left(\mathbf{W}_{\alpha,\beta}[\mu]\right)^q\right](x)&= \int_{0}^{\infty}\left(\frac{1}{t^{N-\alpha \beta}}\int_{B_t(x)}\left(\myint{0}{\infty}\left(\frac{\mu(B_r(y))}{r^{N-\alpha\beta}}\right)^{\frac{1}{\beta-1}}\frac{dr}{r}\right)^qdy\right)^{\frac{1}{\beta-1}}\frac{dt}{t}
	\\
	&\geq c_{_{38}}\int_{0}^{\infty}\left(\frac{1}{t^{N-\alpha\beta}}\int_{B_t(x)}\left(\frac{\mu(B_{2t}(y))}{t^{N-
 			\alpha\beta}}\right)^{\frac{q}{\beta-1}}dy\right)^{\frac{1}{\beta-1}}\frac{dt}{t}
 	\\&\geq c_{_{36}}\int_{0}^{\infty}\left(t^{\frac{\ga\gb(\gb-1) }{q}}\frac{\mu(B_{t}(x))}{t^{N-
 			\alpha\beta}}\right)^{\frac{q}{(\beta-1)^2}}\frac{dt}{t}
 	\\&= c_{_{36}}	\mathbf{W}_{\frac{\alpha\beta(q+\beta-1)}{q+(\beta-1)^2},\frac{(\beta-1)^2}{q}+1}[\mu](x),
 	\end{align*}
 where $c_{_{38}}=c_{_{38}}(\ga,\gb,N,q)>0$, which proves \eqref{wolffes1}. \smallskip

In order to prove \eqref{wolffes1} we recall the following estimate on Wolff potentials \cite{22VHV} 
	\bel{est-W}
 	||\mathbf{W}_{\alpha,\beta}[\omega]||_{L^{\frac{(\beta-1)N}{N-\alpha\beta},\infty}}\leq c_{_{39}}\left(\omega(\mathbb{R}^N)\right)^{\frac{1}{\beta-1}} \qquad\forall\,\omega\in \mathfrak{M}_b^+(\mathbb{R}^N),
 	\ee
 	where $L^{\frac{(p-1)N}{N-\alpha\beta},\infty}$ denotes the weak-$L^{\frac{(p-1)N}{N-\alpha\beta}}$ space.
 	In particular, since $0<q<\frac{N(\beta-1)}{N-\alpha\beta}$, 
 	\bel{est-W'}
 	\myint{B_r(x)}{}\left(\mathbf{W}_{\alpha,\beta}[\omega]\right)^{q}dy\leq c_{_{40}}r^N\left(\frac{\omega(\mathbb{R}^N)}{r^{N-\alpha\beta}}\right)^{\frac{q}{\beta-1}} \qquad\forall x\in\mathbb{R}^N,\,\forall r>0.
 	\ee
 	Applying this inequality to $\omega=\chi_{B_{2r}(x)}\mu$ yields 
 	\begin{align}\label{es2}
 	\myint{B_r(x)}{}\left(\mathbf{W}^r_{\alpha,\beta}[\mu]\right)^{q}dy\leq c_{_{40}}r^N\left(\frac{\mu(B_{2r}(x))}{r^{n-\alpha\beta}}\right)^{\frac{q}{\beta-1}}  \qquad\forall x\in\mathbb{R}^N,\,\forall r>0.
 	\end{align}
 	We claim that
 	\bel{es1}\BA {lll}
 	I:=\myint{0}{\infty}\left(\myfrac{1}{t^{N-\alpha \beta}}\myint{B_t(x)}{}\left(\myint{t}{\infty}\left(\myfrac{\mu(B_r(y))}{r^{N-\alpha\beta}}\right)^{\frac{1}{\beta-1}}\myfrac{dr}{r}\right)^qdy\right)^{\frac{1}{\beta-1}}\myfrac{dt}{t}\\[4mm]
	\phantom{I:}\leq c_{_{37}} \mathbf{W}_{\frac{\alpha\beta(q+\beta-1)}{q+(\beta-1)^2},\frac{(\beta-1)^2}{q}+1}[\mu](x).
 \EA	\ee
 	Since $B_{r}(y)\subset B_{2r}(x)$ for any $y\in B_t(x), r \geq t$,  we have 
$$\BA {lll}
\myint{B_t(x)}{}\left(\myint{t}{\infty}\left(\myfrac{\mu(B_r(y))}{r^{N-\alpha\beta}}\right)^{\frac{1}{\beta-1}}\myfrac{dr}{r}\right)^qdy&\leq \myint{B_t(x)}{}\left(\myint{t}{\infty}\left(\myfrac{\mu(B_{2r}(x))}{r^{N-\alpha\beta}}\right)^{\frac{1}{\beta-1}}\myfrac{dr}{r}\right)^qdy\\[4mm]& \leq c_{_N}t^N\left(\myint{t}{\infty}\left(\myfrac{\mu(B_{2r}(x))}{r^{N-\alpha\beta}}\right)^{\frac{1}{\beta-1}}\myfrac{dr}{r}\right)^q.
\EA$$
 	Hence, 
 	\begin{align*}
I\leq c_{_N}^{\frac{1}{\gb-1}} \int_{0}^{\infty}t^{\frac{\alpha \beta}{\beta-1}}\left(\int_{t}^{\infty}\left(\frac{\mu(B_{2r}(x))}{r^{N-\alpha\beta}}
\right)^{\frac{1}{\beta-1}}\frac{dr}{r}\right)^{\frac{q}{\beta-1}}\frac{dt}{t}.
 	\end{align*}
Using Lemma \ref{hardy}, we infer 
 \begin{align*}
 I&\leq c_{_{37}} \int_{0}^{\infty}r^{\frac{\alpha \beta}{\beta-1}}\left(\frac{\mu(B_{r}(x))}{r^{N-\alpha\beta}}\right)^{\frac{q}{(\beta-1)^2}}\frac{dr}{r}=c_{_{37}}\mathbf{W}_{\frac{\alpha\beta(q+\beta-1)}{q+(\beta-1)^2},\frac{(\beta-1)^2}{q}+1}[\mu](x),
 \end{align*}
which completes the proof. 
 \qeda\medskip
 
The following is a version of Lemma \ref{241220142} for truncated Wolff potentials,
 \begin{lemma} \label{220320152}Let  $1<\beta<N/\alpha$ and $q>0$. If $\delta\in (0,1)$ there holds for any $\gm\in \mathfrak M^+(\BBR^N)$
 	\begin{align}\label{220320153}
 	\mathbf{W}^{\frac{\delta d}{2}}_{\frac{\alpha\beta(q+\beta-1)}{q+(\beta-1)^2},\frac{(\beta-1)^2}{q}+1}[\mu](x)\leq c_{_{42}} \mathbf{W}^{\delta d}_{\alpha,\beta}\left[\left(\mathbf{W}^{\delta d(.)}_{\alpha,\beta}[\mu]\right)^q\right](x)
 	\end{align}
in $\Gw$. Moreover,  if $0< q <\frac{N(\beta-1)}{N-\alpha\beta}$, there holds for any $\gm\in \mathfrak M^+(\BBR^N)$,
 	\begin{align}\label{rev}
 	\mathbf{W}^{R}_{\alpha,\beta}\left[\left(\mathbf{W}^{R}_{\alpha,\beta}[\mu]\right)^q\right](x)\leq c_{_{43}} \mathbf{W}^{4R}_{\frac{\alpha\beta(q+\beta-1)}{q+(\beta-1)^2},\frac{(\beta-1)^2}{q}+1}[\mu](x)
 	\end{align}
	 in $\mathbb{R}^N$.
 \end{lemma}

\nind\proof For any $x\in \Omega$,
 	\begin{align*}\BA {lll}
 	\mathbf{W}^{\delta d(x)}_{\alpha,\beta}\left[\left(\mathbf{W}^{\delta d(.)}_{\alpha,\beta}[\mu](.)\right)^q\right](x)\\\displaystyle
	\phantom{\left[\left(\mathbf{W}^{\delta d(.)}_{\alpha,\beta}[\mu](.)\right)^q\right](x)}
	= \int_{0}^{\delta d(x)}\left(\frac{1}{t^{N-\alpha \beta}}\int_{B_t(x)}\left(\int_{0}^{\delta d(y)}\left(\frac{\mu(B_r(y))}{r^{N-\alpha\beta}}\right)^{\frac{1}{\beta-1}}\frac{dr}{r}\right)^qdy\right)^{\frac{1}{\beta-1}}\frac{dt}{t}.
 	\EA\end{align*}
 	Since $\delta d(y)\geq \frac{7\delta}{8}d(x)$ for all $y\in  B_{\frac{t}{8}}(x)$, provided $0<t<\delta d(x)$,
	\begin{align*}
 	\myint{B_t(x)}{}\left(\int_{0}^{\delta d(y)}\left(\frac{\mu(B_r(y))}{r^{N-\alpha\beta}}\right)^{\frac{1}{\beta-1}}\frac{dr}{r}\right)^qdy &\geq \int_{B_{t/8}(x)}\left(\int_{0}^{\frac{7\delta}{8}d(x)}\left(\frac{\mu(B_r(y))}{r^{N-\alpha\beta}}\right)^{\frac{1}{\beta-1}}\frac{dr}{r}\right)^qdy\\
 	&\geq \int_{B_{t/8}(x)}\left(\int_{0}^{\frac{7t}{8}}\left(\frac{\mu(B_r(y))}{r^{N-\alpha\beta}}\right)^{\frac{1}{\beta-1}}\frac{dr}{r}\right)^qdy
	\\&\geq 
 	c_{_{44}} \int_{B_{t/8}(x)}\left(\frac{\mu(B_{\frac{3t}{4}}(y))}{t^{N-\alpha\beta}}\right)^{\frac{q}{\beta-1}}dy
 	\\&\geq 
 	c_{_{44}} \int_{B_{t/8}(x)}\left(\frac{\mu(B_{\frac{3t}{4}-\frac{t}{8}}(x))}{t^{N-\alpha\beta}}\right)^{\frac{q}{\beta-1}}dy
 	\\&\geq 
 	c_{_{45}} t^{N}\left(\frac{\mu(B_{\frac{t}{2}}(x))}{t^{N-\alpha\beta}}\right)^{\frac{q}{\beta-1}}.
 	\end{align*}
	Hence
	\begin{align*}
	\mathbf{W}^{\delta d(x)}_{\alpha,\beta}\left[\left(\mathbf{W}^{\delta d(.)}_{\alpha,\beta}[\mu](.)\right)^q\right](x)\geq 
	c_{_{46}} \int_{0}^{\delta d(x)}\left(t^{\alpha\beta}\left(\frac{\mu(B_{\frac{t}{2}}(x))}{t^{N-\alpha\beta}}\right)^{\frac{q}{\beta-1}}\right)^{\frac{1}{\beta-1}}\frac{dt}{t},
	\end{align*}
 which implies \eqref{220320153}. \smallskip
 
Because of \eqref{es2}, it is sufficient to prove that there holds
 	\begin{align}\label{es3}
 	\int_{0}^{R}\left(\frac{1}{t^{N-\alpha \beta}}\myint{B_t(x)}{}\left(\int_{t}^{R}\left(\frac{\mu(B_r(y))}{r^{N-\alpha\beta}}\right)^{\frac{1}{\beta-1}}\frac{dr}{r}\right)^qdy\right)^{\frac{1}{\beta-1}}\frac{dt}{t}\leq c_{_{47}} \mathbf{W}^{4R}_{\frac{\alpha\beta(q+\beta-1)}{q+(\beta-1)^2},\frac{(\beta-1)^2}{q}+1}[\mu](x),
\end{align}
in order to obtain \eqref{rev}. 	Since $B_{\rho}(y)\subset B_{2\rho}(x)$ for any $y\in B_r(x), \rho \geq r$,  we have 
 	\begin{align*}
 	\myint{B_t(x)}{}\left(\int_{t}^{R}\left(\frac{\mu(B_r(y))}{r^{N-\alpha\beta}}\right)^{\frac{1}{\beta-1}}\frac{dr}{r}\right)^qdy&\leq \myint{B_t(x)}{}\left(\int_{t}^{R}\left(\frac{\mu(B_{2r}(x))}{r^{N-\alpha\beta}}\right)^{\frac{1}{\beta-1}}\frac{dr}{r}\right)^qdy\\& \leq c_{_N}t^N\left(\int_{t}^{R}\left(\frac{\mu(B_{2r}(x))}{r^{N-\alpha\beta}}\right)^{\frac{1}{\beta-1}}\frac{dr}{r}\right)^q.
\end{align*}
Therefore
 	\begin{align*}
 	&\int_{0}^{R}\left(\frac{1}{t^{N-\alpha \beta}}\myint{B_t(x)}{}\left(\int_{t}^{R}\left(\frac{\mu(B_r(y))}{r^{N-\alpha\beta}}\right)^{\frac{1}{\beta-1}}\frac{dr}{r}\right)^qdy\right)^{\frac{1}{\beta-1}}\frac{dt}{t}\\&\phantom{------------}\leq c_{_N} \int_{0}^{R}\left(t^{\alpha \beta}\left(\int_{t}^{R}\left(\frac{\mu(B_{2r}(x))}{r^{N-\alpha\beta}}\right)^{\frac{1}{\beta-1}}\frac{dr}{r}\right)^q\right)^{\frac{1}{\beta-1}}\frac{dt}{t}.
 	\end{align*}
We infer \eqref{es3} by Lemma \ref{hardy}, which completes the proof.
 \qeda\medskip
 
The next two propositions link Wolff potentials of a measure with Riesz capaciticies (in the case of whole space) and truncated 
Wolff potentials with Bessel capaciticies (in the bounded domain case). Their proof can be found in \cite{22PhVe,22PhVe2} (and \cite{VHV2} with a different method).
\begin{proposition} \label{241020147}Let $1<\beta<N/\alpha$, $q>\beta-1$,  $\nu\in \mathfrak{M}^+(\mathbb{R}^N)$. Then, the following statements are equivalent:\smallskip

\nind (a) The inequality 
		\begin{align}\label{241020141}
		\nu(K)\leq c_{_{48}}\text{Cap}_{\mathbf{I}_{\alpha \beta},\frac{q}{q-\beta+1}}(K)
		\end{align}
		\nind holds for any compact set $K\subset\mathbb{R}^N$, for some $c_{_{48}}>0$.\smallskip

\nind (b) The inequality 
		\begin{align}\label{241020142}
		\int_{\mathbb{R}^N}\left(\mathbf{W}_{\alpha,\beta}[\chi_{B_t(x)}\nu](y)\right)^qdy\leq c_{_{49}} \nu(B_t(x))
		\end{align}
\nind holds for any ball $B_t(x)\subset\mathbb{R}^N$, for some $c_{_{49}}>0$.\smallskip

\nind (c) The inequality 
		\begin{align}\label{241020143}
		\mathbf{W}_{\alpha,\beta}\left[\left(\mathbf{W}_{\alpha,\beta}[\nu]\right)^q\right]\leq c_{_{50}}\mathbf{W}_{\alpha,\beta}[\nu]<\infty ~\text{a.e in }~\mathbb{R}^N
		\end{align}
\nind holds  for some $c_{_{50}}>0$.
\end{proposition}
\begin{proposition}\label{241020148} Let $1<\beta<N/\alpha$, $q>\beta-1$, $R>0$ and $\nu\in \mathfrak{M}_b^+(B_R(x_0))$ for some $x_0\in \mathbb{R}^N$. Then, the following statements are equivalent:\smallskip

\nind (a) The inequality 
		\begin{align}\label{241020144}
		\nu(K)\leq c_{_{51}}\text{Cap}_{\mathbf{G}_{\alpha \beta},\frac{q}{q-\beta+1}}(K)
		\end{align}
\nind holds for any compact set $K\subset\mathbb{R}^N$, for some $c_{_{51}}=c_{_{51}}(R)>0$.\smallskip

\nind (b) The inequality 
		\begin{align}\label{241020145}
		\int_{\mathbb{R}^N}\left(\mathbf{W}^{4R}_{\alpha,\beta}[\chi_{B_t(x)}\nu](y)\right)^qdy\leq c_{_{52}} \nu(B_t(x))
		\end{align}
		\nind holds for any ball $B_t(x)\subset\mathbb{R}^N$, for some $c_{_{52}}=c_{_{52}}(R)>0$.\smallskip

\nind (c) The inequality 
		\begin{align}\label{241020146}
		\mathbf{W}^{4R}_{\alpha,\beta}\left[\left(\mathbf{W}^{4R}_{\alpha,\beta}[\nu]\right)^q\right]\leq c_{_{53}} \mathbf{W}^{4R}_{\alpha,\beta}[\nu]~~\text{ a.e in }~B_{2R}(x_0)
		\end{align}
\nind holds  for some $c_{_{53}}=c_{_{53}}(R)>0$.
\end{proposition}

In the following statement we obtain capacitary estimates on combination of measures. 
  \begin{proposition}\label{261220141} Let $\eta,\mu$ be in $\mathfrak{M}^+(\mathbb{R}^N)$.
  	Assume that $0<q <\frac{N(\beta-1)}{N-\alpha\beta}$ and $qs>(\beta-1)^2$. \smallskip
	
\nind (i) If there holds
  	\begin{align}\label{261220143}
  	\eta(K)+\int_K \left(\mathbf{W}_{\alpha,\beta}[\mu]\right)^sdx\leq  \operatorname{Cap}_{\mathbf{I}_{\frac{\alpha\beta(q+\beta-1)}{q},\frac{qs}{qs-(\beta-1)^2}}}(K),
  	\end{align}
  	for any compact set $K\subset\mathbb{R}^N$, then
  	 \begin{align}\label{250320151}
  	 \mathbf{W}_{\alpha,\beta}\left[\left(\mathbf{W}_{\alpha,\beta}\left[\left(\mathbf{W}_{\alpha,\beta}[\omega]\right)^q\right]\right)^s\right]\leq c_{_{54}} \mathbf{W}_{\alpha,\beta}[\omega]<\infty~~\text{a.e in }~~\mathbb{R}^N,~~
  	 \end{align}
  	 where $\omega=\left(\mathbf{W}_{\alpha,\beta}[\mu]\right)^s+\eta$.\smallskip
	
\nind (ii)  If there holds
  	 \begin{align}\label{261220144}
  	 \eta(K)+\int_K \left(\mathbf{W}^{2R}_{\alpha,\beta}[\mu]\right)^sdx\leq  \operatorname{Cap}_{\mathbf{G}_{\frac{\alpha\beta(q+\beta-1)}{q},\frac{qs}{qs-(\beta-1)^2}}}(K),
  	 \end{align}
  	 for any compact set $K\subset\mathbb{R}^N$, then
  	 \begin{align}\label{260320152}
  	 \mathbf{W}^{2R}_{\alpha,\beta}\left[\left(\mathbf{W}^{2R}_{\alpha,\beta}\left[\left(\mathbf{W}^{2R}_{\alpha,\beta}[\omega]\right)^q\right]\right)^s\right]\leq c_{_{55}} \mathbf{W}^{2R}_{\alpha,\beta}[\omega]<\infty~~\text{a.e in }~~B_{R}(x_0),~~
  	 \end{align}
  	 where $\omega=\chi_{_{B_{R}(x_0)}}\left(\mathbf{W}^{2R}_{\alpha,\beta}[\mu]\right)^s+\chi_{_{B_{R}(x_0)}}\eta$.
  \end{proposition}
  \proof
  	{\it Statement (i)}:	We assume that \eqref{261220143} holds. Put $\omega=\left(\mathbf{W}_{\alpha,\beta}[\mu]\right)^s+\eta$ and apply \eqref{261220143} to $K=\overline{B}_{2\rho}(x)$. Since by homogeneity
  	\begin{align*}
  	\operatorname{Cap}_{\mathbf{I}_{\frac{\alpha\beta(q+\beta-1)}{q},\frac{qs}{qs-(\beta-1)^2}}}(\overline{B}_{2\rho}(x))=\rho^{N-\frac{\alpha\beta(q+\beta-1)s}{qs-(\beta-1)^2}}\operatorname{Cap}_{\mathbf{I}_{\frac{\alpha\beta(q+\beta-1)}{q},\frac{qs}{qs-(\beta-1)^2}}}\!\!\!\!\!\!(\overline{B}_{2}(0)),
  	\end{align*} 
	we deduce from \eqref{261220143}
  		\begin{align*}
  		\omega(B_\rho(x))\leq c_{_{55}} \rho^{N-\frac{\alpha\beta(q+\beta-1)s}{qs-(\beta-1)^2}}\qquad\forall~~\rho>0,
  		\end{align*} 
which is equivalent to 
  		\begin{align}\label{241220144}
  		\rho^{\frac{\alpha\beta}{\beta-1}}\left(\frac{\omega(B_{\rho}(x))}{\rho^{N-\frac{\alpha\beta(q+\beta-1)}{q}}}\right)^{\frac{qs}{(\beta-1)^3}}
		\leq c_{_{56}} \left(\frac{\omega(B_{\rho
  			}(x))}{\rho^{N-\alpha\beta}}\right)^{\frac{1}{\beta-1}}\qquad\forall~~\rho>0.
  		\end{align}
  		We apply Proposition \ref{241020147} to $\nu=\omega$ with $(\alpha,\beta,q)=\left(\frac{\alpha\beta(q+\beta-1)}{q+(\beta-1)^2},\frac{(\beta-1)^2}{q}+1,s\right)$, \eqref{261220143}  implies 
  		\begin{align}\label{241220145}
  		\int_{\mathbb{R}^N}\left(\mathbf{W}_{\frac{\alpha\beta(q+\beta-1)}{q+(\beta-1)^2},\frac{(\beta-1)^2}{q}+1}[\chi_{_{B_{t}(x)}}\omega]\right)^sdy\leq c_{_{57}}\omega(B_t(x)).
  		\end{align}
  		By Lemma \ref{241220142},  
  		\eqref{250320151} is equivalent to  
  		\begin{align}\label{241220146}
  		\mathbf{W}_{\alpha,\beta}\left[\left(\mathbf{W}_{\frac{\alpha\beta(q+\beta-1)}{q+(\beta-1)^2},\frac{(\beta-1)^2}{q}+1}[\omega]\right)^s\right]\leq c_{_{58}} \mathbf{W}_{\alpha,\beta}[\omega]<\infty~~\text{a.e}~~\mathbb{R}^N.
  		\end{align}	
  		Therefore, it is enough to show that \eqref{241220144} and \eqref{241220145} imply \eqref{241220146}. In fact, since for $t>0$
  		\begin{align*}
  		\myint{B_t(x)}{}\left(\mathbf{W}^t_{\frac{\alpha\beta(q+\beta-1)}{q+(\beta-1)^2},\frac{(\beta-1)^2}{q}+1}[\omega](y)\right)^sdy= \myint{B_t(x)}{}\left(\mathbf{W}^t_{\frac{\alpha\beta(q+\beta-1)}{q+(\beta-1)^2},\frac{(\beta-1)^2}{q}+1}[\chi_{_{B_{2t}(x)}}\omega](y)\right)^sdy,
  		\end{align*}
  		we apply  \eqref{241220145} and obtain 
  		\begin{align*}
  		\myint{B_t(x)}{}\left(\mathbf{W}^t_{\frac{\alpha\beta(q+\beta-1)}{q+(\beta-1)^2},\frac{(\beta-1)^2}{q}+1}[\omega](y)\right)^sdy\leq c_{_{57}} \omega(B_{2t}(x)).
  		\end{align*}
  		So, it is enough to show that 
  		\begin{align}\label{241220147}
  		I:=\int_{0}^{\infty}\left(\frac{1}{t^{N-\alpha \beta}}\myint{B_t(x)}{}\left(\int_{t}^{\infty}\left(\frac{\omega(B_r(y))}{r^{N-\frac{\alpha\beta(q+\beta-1)}{q}}}\right)^{\frac{q}{(\beta-1)^2}}\frac{dr}{r}\right)^sdy\right)^{\frac{1}{\beta-1}}\frac{dt}{t}\leq c_{_{58}} \mathbf{W}_{\alpha,\beta}[\omega](x).
  		\end{align}
  		Since $B_{r}(y)\subset B_{2r}(x)$ for any $y\in B_t(x), r \geq t$,  we have 
  		\begin{align*}
  		I&\leq c_{_N}\int_{0}^{\infty}\left(t^{\alpha \beta}\left(\int_{t}^{\infty}\left(\frac{\omega(B_{2r}(x))}{r^{N-\frac{\alpha\beta(q+\beta-1)}{q}}}\right)^{\frac{q}{(\beta-1)^2}}\frac{dr}{r}\right)^s\right)^{\frac{1}{\beta-1}}\frac{dt}{t}
  		\\& =c_{_N}\int_{0}^{\infty}t^{\frac{\alpha\beta}{\beta-1}}\left(\int_{t}^{\infty}\left(\frac{\omega(B_{2r}(x))}{r^{N-\frac{\alpha\beta(q+\beta-1)}{q}}}\right)^{\frac{q}{(\beta-1)^2}}\frac{dr}{r}\right)^{\frac{s}{\beta-1}}\frac{dt}{t}.
  		\end{align*}
  		It follows from  Lemma \ref{hardy} and \eqref{241220144} that 
  			\begin{align*}
  		I\leq c_{_{59}}\int_{0}^{\infty}t^{\frac{\alpha\beta}{\beta-1}}\left(\frac{\omega(B_{2t}(x))}{t^{N-\frac{\alpha\beta(q+\beta-1)}{q}}}\right)^{\frac{qs}{(\beta-1)^3}}\frac{dt}{t}\leq c_{_{56}}c_{_{59}}\int_{0}^{\infty}\left(\frac{\omega(B_{2t}(x))}{t^{N-\alpha\beta}}\right)^{\frac{1}{\beta-1}}\frac{dt}{t},
  		\end{align*}
which is \eqref{241220147}.\smallskip

\nind   	{\it Statement (ii)}:	We assume that \eqref{261220144} holds. Put $d\omega=\chi_{_\Omega}\left(\mathbf{W}_{\alpha,\beta}[\mu]\right)^s+\chi_{_\Omega}\eta$, then 
  			\begin{align*}
  			\omega(B_\rho(x))\leq c_{_{60}} \rho^{N-\frac{\alpha\beta(q+\beta-1)s}{qs-(\beta-1)^2}}\qquad\forall~0<\rho<2R.
  			\end{align*} 
  			As in the proof of statement (i),  the above inequality is equivalent to 
  			\begin{align}\label{2203201513}
  			\rho^{\frac{\alpha\beta}{\beta-1}}\left(\frac{\omega(B_{\rho}(x))}{\rho^{N-\frac{\alpha\beta(q+\beta-1)}{q}}}\right)^{\frac{qs}{(\beta-1)^3}}\leq c_{_{61}} \left(\frac{\omega(B_{\rho
  				}(x))}{\rho^{N-\alpha\beta}}\right)^{\frac{1}{\beta-1}}\qquad\forall~0<\rho<2R.
  			\end{align}
  				Applying Proposition \ref{241020148} with $\nu=\omega$ and $(\alpha,\beta,q)=\left(\frac{\alpha\beta(q+\beta-1)}{q+(\beta-1)^2},\frac{(\beta-1)^2}{q}+1,s\right)$, 
  			\begin{align}\label{2203201514}
  			\int_{\mathbb{R}^N}\left(\mathbf{W}^{4R}_{\frac{\alpha\beta(q+\beta-1)}{q+(\beta-1)^2},\frac{(\beta-1)^2}{q}+1}[\chi_{_{B_{t}(x)}}\omega]\right)^sdy\leq c_{_{62}}\omega(B_t(x)).
  			\end{align}By Lemma \ref{220320152}, \eqref{260320152} is equivalent to 
  			\begin{align}\label{2203201515}
  			\mathbf{W}^{4R}_{\alpha,\beta}\left[\left(\mathbf{W}^{4R}_{\frac{\alpha\beta(q+\beta-1)}{q+(\beta-1)^2},\frac{(\beta-1)^2}{q}+1}[\omega]\right)^s\right]\leq c_{_{63}} \mathbf{W}^{4R}_{\alpha,\beta}[\omega]\qquad\text{a.e in}~B_{R}(x_0).
  			\end{align}
  			Therefore, it is sufficient to prove that \eqref{2203201513} and \eqref{2203201514} imply \eqref{2203201515}. Actually, since 
  			\begin{align*}
  			\myint{B_t(x)}{}\left(\mathbf{W}^t_{\frac{\alpha\beta(q+\beta-1)}{q+(\beta-1)^2},\frac{(\beta-1)^2}{q}+1}[\omega](y)\right)^sdy= \myint{B_t(x)}{}\left(\mathbf{W}^t_{\frac{\alpha\beta(q+\beta-1)}{q+(\beta-1)^2},\frac{(\beta-1)^2}{q}+1}[\chi_{B_{2t}(x)}\omega](y)\right)^sdy
  			\end{align*}
  			for all $0<t<4R$, thus applying  \eqref{2203201514}, we obtain 
  			\begin{align*}
  			\myint{B_t(x)}{}\left(\mathbf{W}^t_{\frac{\alpha\beta(q+\beta-1)}{q+(\beta-1)^2},\frac{(\beta-1)^2}{q}+1}[\omega](y)\right)^sdy\leq c_{_{64}} \omega(B_{2t}(x)).
  			\end{align*}
  			So, it is sufficient to show that for any $x\in B_{R}(x_0)$
  			\begin{align}\label{2203201516}
  				II:=\int_{0}^{4R}\left(\frac{1}{t^{N-\alpha \beta}}\myint{B_t(x)}{}\left(\int_{t}^{4R}\left(\frac{\omega(B_r(y))}{r^{N-\frac{\alpha\beta(q+\beta-1)}{q}}}\right)^{\frac{q}{(\beta-1)^2}}\frac{dr}{r}\right)^sdy\right)^{\frac{1}{\beta-1}}\frac{dt}{t}\leq c_{_{65}} \mathbf{W}^{4R}_{\alpha,\beta}[\omega](x).
  			\end{align}
  		Since $B_{r}(y)\subset B_{2r}(x)$ for any $y\in B_t(x)$ with $r \geq t$,  we have 
  		\begin{align*}
  		II&\leq c_{_N}\int_{0}^{4R}t^{\frac{\alpha\beta}{\beta-1}}\left(\int_{t}^{4R}\left(\frac{\omega(B_{2r}(x))}{r^{n-\frac{\alpha\beta(q+\beta-1)}{q}}}\right)^{\frac{q}{(\beta-1)^2}}\frac{dr}{r}\right)^{\frac{s}{\beta-1}}\frac{dt}{t}.
  		\end{align*}
  		Combining this with Lemma \ref{hardy} and \eqref{2203201513} yields 
  		\begin{align*}
  		II\leq c_{_{66}}\mathbf{W}^{16R}_{\alpha,\beta}[\omega](x).
  		\end{align*}
  			Therefore, \eqref{2203201515} follows since $\mathbf{W}^{16R}_{\alpha,\beta}[\omega]\leq c_{_{67}} \mathbf{W}^{4R}_{\alpha,\beta}[\omega]$ in $B_{R}(x_0)$. 
  \qeda\medskip
  
\begin{proposition} \label{pro1} Let $\eta,\mu$ be in $\mathfrak{M}^+(\mathbb{R}^N)$.
	Assume that $0<q <\frac{N(\beta-1)}{N-\alpha\beta}$ and $qs>(\beta-1)^2$. Let $(u_m,v_m)$ be nonnegative measurable funtions in $\mathbb{R}^N$ verifying, for all $m\geq 0$,
	\begin{align}\nonumber
	u_{m+1}\leq c^* \mathbf{W}_{\alpha,\beta}[v_m^q+\mu], ~~~~
	v_{m+1}\leq c^* \mathbf{W}_{\alpha,\beta}[u_m^s+\eta]\qquad\text{a.e. in }~~\mathbb{R}^N,
	\end{align}
	 for some $c^*>0$ and $(u_0,v_0)=0$. 
	Then, there exists a constant $M^*>0$ depending only on $N,\alpha,\beta,q,s, c^*$ such that if the measure $d\omega=\left(\mathbf{W}_{\alpha,\beta}[\mu]\right)^sdx+d\eta$ satisfies
		\begin{align}\label{261220143*}
	\omega(K)\leq M^* \operatorname{Cap}_{\mathbf{I}_{\frac{\alpha\beta(q+\beta-1)}{q}},\frac{qs}{qs-(\beta-1)^2}}(K),
	\end{align}
	for any compact set $K\subset\mathbb{R}^N$, then 
	\begin{align}\label{es4"}
	v_m\leq c_{_{69}} \mathbf{W}_{\alpha,\beta}[\omega],~~~u_m\leq c_{_{70}} \mathbf{W}_{\alpha,\beta}[\left(\mathbf{W}_{\alpha,\beta}[\omega]\right)^q]+c_{_{68}} \mathbf{W}_{\alpha,\beta}[\mu]~~\forall~~m\geq 0,
	\end{align}
	for some constants $c_{_{68}},c_{_{69}},c_{_{70}}$ depending only on $N,\alpha,\beta,q,s$ and $c^*$.
\end{proposition}
\proof By Proposition \ref{261220141}, \eqref{261220143*} implies 
	\begin{align}\label{250320151'}
	\mathbf{W}_{\alpha,\beta}\left[\left(\mathbf{W}_{\alpha,\beta}\left[\left(\mathbf{W}_{\alpha,\beta}[\omega]\right)^q\right]\right)^s\right]\leq c_{_{71}} M^{\frac{qs}{(\beta-1)^3}} \mathbf{W}_{\alpha,\beta}[\omega]<\infty\qquad\text{a.e in}~~\mathbb{R}^N.
	\end{align}
	We set
	\begin{align*}
	&c_{_{68}}= c^*2^{\frac{1}{\beta-1}}, \\
	& c_{_{69}}= c^*2^{1+\frac{1}{\beta-1}}(c_{_{68}}^s2^{s-1}+1)^{\frac{1}{\beta-1}},\\ &
	c_{_{70}}=c^*2^{\frac{1}{\beta-1}} c_{_{69}}^{\frac{q}{\beta-1}},
	\end{align*}
	and choose $M^*>0$ such that 
	\begin{align*}
	 c^* 2^{\frac{1}{\beta-1}}\left(c_{_{70}}^{s}2^{s-1}\right)^{\frac{1}{\beta-1}}c_{_{71}} M^{*\,\frac{qs}{(\beta-1)^3}}=\frac{c_{_{69}}}{2}.
	\end{align*}
	We claim that 
	\begin{align}\label{es4}
	v_m\leq c_{_{69}} \mathbf{W}_{\alpha,\beta}[\omega],~~~u_m\leq c_{_{70}} \mathbf{W}_{\alpha,\beta}[\left(\mathbf{W}_{\alpha,\beta}[\omega]\right)^q]+c_{_{68}} \mathbf{W}_{\alpha,\beta}[\mu]\qquad\forall~~m\geq 0.
	\end{align}
	Clearly, by definition of $c_{_{68}},c_{_{69}}$ and $c_{_{70}}$, we have  \eqref{es4} for $m=0,1$. Next we assume that \eqref{es4} holds for all integer $m\leq l$ for some $l\in\BBN_+^{*}$, then 
	\begin{align*}
		u_{l+1}&\leq c^* \mathbf{W}_{\alpha,\beta}[v_l^q+\mu]\\&\leq c^*2^{\frac{1}{\beta-1}} c_{_{69}}^{\frac{q}{\beta-1}}\mathbf{W}_{\alpha,\beta}[\left(\mathbf{W}_{\alpha,\beta}[\omega]\right)^q] +c^*2^{\frac{1}{\beta-1}} \mathbf{W}_{\alpha,\beta}[\mu]\\&= c_{_{70}} \mathbf{W}_{\alpha,\beta}[\left(\mathbf{W}_{\alpha,\beta}[\omega]\right)^q]+c_{_{68}} \mathbf{W}_{\alpha,\beta}[\mu],
	\end{align*}
	and 
	\begin{align*}
		v_{l+1}& \leq c^* \mathbf{W}_{\alpha,\beta}[\left(c_{_{70}} \mathbf{W}_{\alpha,\beta}[\left(\mathbf{W}_{\alpha,\beta}[\omega]\right)^q]+c_{_{68}} \mathbf{W}_{\alpha,\beta}[\mu]\right)^s+\eta]
		\\& \leq c^* \mathbf{W}_{\alpha,\beta}[c_{_{70}}^{s}2^{s-1}\left( \mathbf{W}_{\alpha,\beta}[\left(\mathbf{W}_{\alpha,\beta}[\omega]\right)^q]\right)^s+c_{_{68}}^s2^{s-1} (\mathbf{W}_{\alpha,\beta}[\mu])^s+\eta]
		\\& \leq c^* 2^{\frac{1}{\beta-1}}\left(c_{_{70}}^{s}2^{s-1}\right)^{\frac{1}{\beta-1}}\mathbf{W}_{\alpha,\beta}[\left( \mathbf{W}_{\alpha,\beta}[\left(\mathbf{W}_{\alpha,\beta}[\omega]\right)^q]\right)^s]\\
		&~~~~+c^*2^{\frac{1}{\beta-1}}(c_{_{68}}^s2^{s-1}+1)^{\frac{1}{\beta-1}}  \mathbf{W}_{\alpha,\beta}[(\mathbf{W}_{\alpha,\beta}[\mu])^s+\eta]
		\\& \leq c^* 2^{\frac{1}{\beta-1}}\left(c_{_{70}}^{s}2^{s-1}\right)^{\frac{1}{\beta-1}}c_{_{71}} M^{*\,\frac{qs}{(\beta-1)^3}} \mathbf{W}_{\alpha,\beta}[\omega]+c^*2^{\frac{1}{\beta-1}}(c_{_{68}}^s2^{s-1}+1)^{\frac{1}{\beta-1}}  \mathbf{W}_{\alpha,\beta}[\omega]
		\\& = \frac{c_{_{69}}}{2} \mathbf{W}_{\alpha,\beta}[\omega]+\frac{c_{_{69}}}{2} \mathbf{W}_{\alpha,\beta}[\omega]\\&= c_{_{69}}\mathbf{W}_{\alpha,\beta}[\omega].
	\end{align*}
	Thus, \eqref{es4} holds true for $m=l+1$. Hence, \eqref{es4} is valid for all $l\geq 0.$
\qeda
\medskip

The next result  is an adaptation of Proposition \ref{pro1} to truncated Wolff potentials.
\begin{proposition} \label{pro2} Let $\eta,\mu$ be in $\mathfrak{M}_b^+(B_R(x_0))$.
	Assume that $0<q <\frac{N(\beta-1)}{N-\alpha\beta}$ and $qs>(\beta-1)^2$. Let $(u_m,v_m)$ be nonnegative measurable funtions in $\mathbb{R}^N$ such that for all $m\geq 0$
	\begin{align}\nonumber
	u_{m+1}\leq c_* \mathbf{W}_{\alpha,\beta}^R[\chi_{_{B_{R}(x_0)}}v_m^q+\mu], ~~~~
	v_{m+1}\leq c_* \mathbf{W}_{\alpha,\beta}^R[\chi_{_{B_{R}(x_0)}}u_m^s+\eta]\quad\text{a.e. in }\; B_R(x_0),
	\end{align}
		and $(u_0,v_0)=0$. If we set $d\omega=\left(\mathbf{W}_{\alpha,\beta}^{2R}[\mu]\right)^sdx+d\eta$, there exists a constant $M_*>0$ depending only on $N,\alpha,\beta,q,s, R$ and $ c_*$ such that if 
	\begin{align}\label{261220143**}
	\omega(K)\leq M_* \operatorname{Cap}_{\mathbf{G}_{\frac{\alpha\beta(q+\beta-1)}{q}},\frac{qs}{qs-(\beta-1)^2}}(K),
	\end{align}
	for any compact set $K\subset\mathbb{R}^N$, then 
\begin{align}\label{es4*}
v_m\leq c_{_{73}} \mathbf{W}_{\alpha,\beta}^{2R}[\omega],~~~u_m\leq c_{_{74}}  \mathbf{W}_{\alpha,\beta}^{2R}[\left(\mathbf{W}^{2R}_{\alpha,\beta}[\omega]\right)^q]+c_{_{72}}  \mathbf{W}_{\alpha,\beta}^{2R}[\mu]\qquad\forall~k\geq 0
\end{align}
in $B_R(x_0)$ for some constants $c_{_{72}} ,c_{_{73}} ,c_{_{74}} $ depending only on $N,\alpha,\beta,q,s,R$ and  $c_*$.
\end{proposition}
\proof The proof is similar to the one of Proposition \ref{pro1} and we omit the details.
\qeda\medskip

  \begin{proposition}\label{260320151} Let $1<\gb<N/\ga$ and $q,s>0$ such that $qs>(\gb-1)^2$.\smallskip
  
\nind (i) Assume that $\eta$ and $\mu$ belong to $\mathfrak{M}_b^+(\BBR^N)$ and $(u,v)$ are  nonnegative measurable functions satisfying 
  	\begin{equation}\label{2203201571}
  	\BA
  	{llll}%
(i)\qquad \qquad \qquad  	&\mathbf{W}_{\alpha,\beta}[v^q]+\mathbf{W}_{\alpha,\beta}[\mu]\leq c_{_{75}} u,\qquad \\[2mm]
  (ii)\qquad  	&\mathbf{W}_{\alpha,\beta}[u^s]+\mathbf{W}_{\alpha,\beta}[\eta]\leq c_{_{75}} v \qquad\text{a.e. in }\;\mathbb{R}^N,\qquad                                                             
  	\EA
  	\end{equation}	
for some $c_{_{75}}>0$. Then  there exists a constant $c_{_{76}}>0$ depending only on $N,\alpha,\beta,q,s$ and $c_{_{75}}$ such that  
\begin{align}\label{2203201581}
  	\eta(K)+\int_K \left(\mathbf{W}_{\alpha,\beta}[\mu](x)\right)^sdx\leq c_{_{76}} \operatorname{Cap}_{\mathbf{I}_{\frac{\alpha\beta(q+\beta-1)}{q}},\frac{qs}{qs-(\beta-1)^2}}(K),
  	\end{align}
  	for any compact set $K\subset\BBR^N$.\smallskip
  
\nind (ii) Assume that $\eta$ and $\mu$ belong to $\mathfrak{M}_b^+(\Gw)$ and  $(u,v)$ are  nonnegative functions satisfying 
  	\begin{equation}\label{220320157}
\BA{llll}%
(i)\qquad \qquad \qquad  &\mathbf{W}^{\delta d(.)}_{\alpha,\beta}[v^q]+\mathbf{W}^{\delta d}_{\alpha,\beta}[\mu] \leq c_{_{77}} u,\\[2mm]
(ii)\qquad \qquad \qquad  &\mathbf{W}^{\delta d(.)}_{\alpha,\beta}[u^s]+\mathbf{W}^{\delta d}_{\alpha,\beta}[\eta]\leq c_{_{77}} v \qquad\text{a.e. in }\;\Gw,\qquad                                                                
	\EA  
  	\end{equation} 
for some $c_{_{77}}>0$.	Then for any $\Omega'\subset\subset \Omega$, there exists a constant $c_{_{78}}>0$ depending only on $n,\alpha,\beta,q,s,c_{_{77}}$ and $dist(\Omega',\partial\Omega) $ such that  
	\begin{align}\label{220320158}
  	\eta(K)+\int_K \left(\mathbf{W}^{\delta d(x)}_{\alpha,\beta}[\mu](x)\right)^sdx\leq c_{_{78}} \operatorname{Cap}_{\mathbf{G}_{\frac{\alpha\beta(q+\beta-1)}{q}},\frac{qs}{qs-(\beta-1)^2}}(K),
  	\end{align}
  	for any compact set $K\subset\Omega'$.
  \end{proposition}
  \proof (i): Set $\omega=u^s+\eta$, then 
  	\begin{align*}
  	\omega\geq u^s \geq \left(\mathbf{W}_{\alpha,\beta}[v^q]\right)^s\geq c_{_{79}}\left(\mathbf{W}_{\alpha,\beta}\left[\left(\mathbf{W}_{\alpha,\beta}[\omega]\right)^q\right]\right)^s.
  	\end{align*}
  	By \eqref{wolffes1} in Lemma \ref{241220142}, we get 
  	\begin{align*}
  	\omega\geq c_{_{80}}\left(\mathbf{W}_{\frac{\alpha\beta(q+\beta-1)}{q+(\beta-1)^2},\frac{(\beta-1)^2}{q}+1}[\omega]\right)^s,
  	\end{align*}
  	which implies 
  \begin{align*}
  \int_{\mathbb{R}^N}\left(\mathbf{W}_{\frac{\alpha\beta(q+\beta-1)}{q+(\beta-1)^2},\frac{(\beta-1)^2}{q}+1}[\chi_{_{B_{t}(x)}}\omega]\right)^sdy\leq c_{_{81}}\omega(B_t(x))\qquad\forall \; x\in\mathbb{R}^N,\;\forall\;t>0.
  \end{align*}
 Applying Proposition \ref{241020147} to $\mu=\omega$ with $(\alpha,\beta,q)=\left(\frac{\alpha\beta(q+\beta-1)}{q+(\beta-1)^2},\frac{(\beta-1)^2}{q}+1,s\right)$,  we get \eqref{2203201581}.\smallskip
 
\nind (ii) We define $\gw$ as above and we have 
\begin{align*}
\omega\geq u^s \geq \left(\mathbf{W}^{\delta d}_{\alpha,\beta}[v^q]\right)^s\geq c_{_{82}}\left(\mathbf{W}^{\delta d}_{\alpha,\beta}\left[\left(\mathbf{W}^{\delta d(.)}_{\alpha,\beta}[\omega]\right)^q\right]\right)^s\qquad \text{a.e. in } \Omega,
\end{align*}
which leads to
\begin{align*}
\omega\geq c_{_{83}}\left(\mathbf{W}^{\frac{\delta}{2}d}_{\frac{\alpha\beta(q+\beta-1)}{q+(\beta-1)^2},\frac{(\beta-1)^2}{q}+1}[\omega]\right)^s\qquad \text{a.e. in } \Omega,
\end{align*}
by inequality \eqref{220320153} in Lemma \ref{220320152}. Let $M_\omega$ denote the centered Hardy-Littlewood maximal function which is defined for any $f\in L_{loc}^1(\mathbb{R}^N,d\omega)$  by 
\[{M_\omega }f(x) = \sup _{t > 0} \frac{1}{\omega (B_t(x))}\myint{B_t(x)}{}|f|d\omega  .\]
Let $K\subset\Omega$ be compact. Set $r_K=dist(K,\partial \Omega)$ and $\Omega_K=\{x\in\Omega: d(x,K)<r_K/2\}$.
Then, for any Borel set $E\subset K$,  
\[c_{_{84}}\int_\Omega  {{{\left( {{M_\omega }{\chi _{_E}}} \right)}^{\frac{{sq }}{{(\beta-1)^2}}}}\left(\mathbf{W}^{\frac{\delta}{2}d(x)}_{\frac{\alpha\beta(q+\beta-1)}{q+(\beta-1)^2},\frac{(\beta-1)^2}{q}+1}[\omega]\right)^sdx}  \le  \int_\Omega  {{{\left( {{M_\omega }{\chi _{_E}}} \right)}^{\frac{{sq }}{{(\beta-1)^2}}}}d\omega }. \]
Since $M_\omega$ is a bounded linear map on $L^p(\mathbb{R}^N,d\omega)$ for any $p>1$ and 
$$\BA {lll}\left( M_\omega \chi _{_E} \right)^{\frac{sq }{(\beta-1)^2}}\left(\mathbf{W}^{\frac{\delta}{2}d(x)}_{\frac{\alpha\beta(q+\beta-1)}{q+(\beta-1)^2},\frac{(\beta-1)^2}{q}+1}[\omega]\right)^s\geq \myint{0}{\frac{\delta}{2}d(x)}\left(\myfrac{\gw (B_t(x)\cap E)}{\gw(B_t(x))}
\myfrac{\gw(B_t(x)}{t^{N-\frac{\alpha\beta(q+\beta-1)}{q}}}\right)^{\frac{sq}{(\gb-1)^2}}\myfrac{dt}{t},
\EA$$
we obtain
\[\int_\Omega  {\left(\mathbf{W}^{\frac{\delta}{2}d(x)}_{\frac{\alpha\beta(q+\beta-1)}{q+(\beta-1)^2},\frac{(\beta-1)^2}{q}+1}[\omega_E]\right)^sdx}  \le  c_{_{85}}{\omega }(E),\]
where $\omega_E=\chi_{_E}\omega$. 
Note that if $x\in \Omega$ and $d(x)\le r_K/8$, then $B_t(x)\subset \Omega\backslash \Omega_K$ for all $t\in (0,\frac{\delta d(x)}{2})$;
indeed, for all $y\in B_t(x)$\\
\[d(y,\partial \Omega ) \le d(x,\partial \Omega ) + |x - y| < (1 + \delta )d(x,\partial \Omega ) < \frac{1}{4}{r_K},\]
thus
\[d(y,K) \geq d(K,\partial \Omega ) - d(y,\partial \Omega ) > \frac{3}{4}{r_K} > \frac{1}{2}{r_K},\]
which implies $y\notin \Omega_K$. 
We deduce that  
\[ {\bf W}_{\frac{\alpha\beta(q+\beta-1)}{q+(\beta-1)^2},\frac{(\beta-1)^2}{q}+1}^{\frac{\delta}{2} d(x,\partial \Omega )}[{\omega _E}](x) \geq  {\bf W}_{\frac{\alpha\beta(q+\beta-1)}{q+(\beta-1)^2},\frac{(\beta-1)^2}{q}+1}^{\frac{\delta }{16}{r_K}}[{\omega _E}](x)\qquad\forall  x\in \Omega,\]
and 
\[ {\bf W}_{\frac{\alpha\beta(q+\beta-1)}{q+(\beta-1)^2},\frac{(\beta-1)^2}{q}+1}^{\frac{\delta }{16}{r_K}}[{\omega _E}](x) = 0\qquad\forall  x\in \Omega^c.\]
Hence we obtain 
\begin{equation}\label{2hvine4} \int_{\mathbb{R}^N} {\left(\mathbf{W}^{\frac{\delta}{16}r_K}_{\frac{\alpha\beta(q+\beta-1)}{q+(\beta-1)^2},\frac{(\beta-1)^2}{q}+1}[\omega_E]\right)^sdx}  \le  c_{_{85}}{\omega }(E)\qquad\forall  E\subset K,\, E\text{ Borel}.
\end{equation}
Applying Proposition \ref{241020148}  with $\mu=\chi_{_{K\cap B_{2^{-6}\delta r_{_K}}\!\!\!\!\!\!(x) }}\omega$  we get \eqref{220320158}, which completes the proof. $\phantom{--------}$
  \qeda
  

  \section{Quasilinear Dirichlet problems}
  Let $\Omega$ be a bounded domain in $\mathbb{R}^N$.   If $\mu\in\mathfrak{M}_b(\Omega)$, we denote by $\mu^+$ and $\mu^-$ respectively its positive and negative parts in the Jordan decomposition. We denote by $\mathfrak{M}_0(\Omega)$ the space of measures in $\Omega$ which are absolutely continuous with respect to the $c^{\Omega}_{1,p}$-capacity defined on a compact set $K\subset\Omega$ by
  \begin{equation*}
  c^{\Omega}_{1,p}(K)=\inf\left\{\int_{\Omega}{}|{\nabla \varphi}|^pdx:\varphi\geq \chi_{_K},\varphi\in C^\infty_c(\Omega)\right\}.
  \end{equation*}
  We also denote $\mathfrak{M}_s(\Omega)$ the space of measures in $\Omega$ with support on a set of zero $c^{\Omega}_{1,p}$-capacity. Classically, any $\mu\in\mathfrak{M}_b(\Omega)$ can be written in a unique way under the form $\mu=\mu_0+\mu_s$ where $\mu_0\in \mathfrak{M}_0(\Omega)\cap \mathfrak{M}_b(\Omega)$ and $\mu_s\in \mathfrak{M}_s(\Omega)\cap \mathfrak{M}_b(\Omega)$.
  It is well known  that any  $\mu_0\in \mathfrak{M}_0(\Omega)\cap\mathfrak{M}_b(\Omega)$ can be written under the form $\mu_0=f-div ~g$ where $f\in L^1(\Omega)$ and $g\in L^{p'}(\Omega,\mathbb{R}^N)$.
  
  For $k>0$ and $s\in\mathbb{R}$ we set $T_k(s)=\max\{\min\{s,k\},-k\}$. If $u$ is a measurable function defined  in $\Omega$, finite a.e. and such that $T_k(u)\in W^{1,p}_{loc}(\Omega)$ for any $k>0$, there exists a measurable function $v:\Omega\to \mathbb{R}^N$ such that $\nabla T_k(u)=\chi_{_{\{|u|\leq k\}}}v$ 
  a.e. in $\Omega$ and for all $k>0$. We define the gradient a.e. $\nabla u$ of $u$ by $v=\nabla u$. We recall the definition of a renormalized solution given in \cite{22DMOP}.
  
  \begin{definition} Let $\mu=\mu_0+\mu_s\in\mathfrak{M}_b(\Omega)$. A measurable  function $u$ defined in $\Omega$ and finite a.e. is called a renormalized solution of 
  	
  	\begin{equation}
  	\label{2hvpro1}\begin{array}{ll}
  	- {\Delta _p}u = \mu \qquad&\;in\;\Omega  \\ 
  	\phantom{ - {\Delta _p}}u = 0~~~&\;on\;\partial \Omega,  \\ 
  	\end{array}
  	\end{equation}
  	if $T_k(u)\in W^{1,p}_0(\Omega)$ for any $k>0$, $|{\nabla u}|^{p-1}\in L^r(\Omega)$ for any $0<r<\frac{N}{N-1}$, and $u$ has the property that for any $k>0$ there exist $\lambda_k^+$ and $\lambda_k^-$ belonging to $\mathfrak{M}_{b}^+\cap\mathfrak{M}_0(\Omega)$, respectively concentrated on the sets $u=k$ and $u=-k$, with the property that 
  	$\mu_k^+\rightharpoonup\mu_s^+$, $\mu_k^-\rightharpoonup\lambda_s^-$ in the narrow topology of measures and such that
  	\[
  	\int_{\{|u|<k\}}\left\vert \nabla u\right\vert ^{p-2}\nabla u.\nabla\varphi
  	dx=\int_{\{|u|<k\}}{\varphi d}{\mu_{0}}+\int_{\Omega}\varphi d\lambda_{k}%
  	^{+}-\int_{\Omega}\varphi d\lambda_{k}^{-},%
  	\]
  	for every $\varphi\in W^{1,p}_0(\Omega)\cap L^{\infty}(\Omega)$.
  	
  \end{definition}
  
  \begin{remark}
  	\label{2hvR1} We recall that if $u$ is a renormalized solution to problem \eqref{2hvpro1}, then $\frac{|\nabla u|^p}{(|u|+1)^r}\in L^1(\Omega)$ for all $r>1$.  Furthermore, $u\geq 0$ $a.e.$ in $\Omega$ if $\mu\in\mathfrak{M}_{b}^+(\Omega)$.
  \end{remark}
  
  The following general stability result has been proved in \cite[Th 4.1]{22DMOP}.
  
  \begin{theorem}
  	\label{2hvP3} Let $\mu=\mu_{0}+\mu_{s}^{+}-\mu_{s}^{-},$ with $\mu_{0}%
  	=F-\operatorname{div}g\in\mathfrak{M}_{0}(\Omega)$ and  $\mu_{s}^{+}$, $\mu ^{-}_{s}$
  	belonging to $\mathfrak{M}_{s}^{+}(\Omega).$ Let
  	$ \mu_{n}=F_{n}-\operatorname{div}g_{n}+\rho_{n}-\eta_{n}$ with
  	$F_{n}\in L^{1}(\Omega)$, $g_{n} \in(L^{p^{\prime}}(\Omega))^{N}$ and $\rho_{n}$, $\eta
  	_{n}$ belonging to $\mathfrak{M}_{b}^{+}(\Omega)$.
  	Assume that $\{F_{n}\}$ converges to $F$ weakly in $L^{1}(\Omega)$, $\{g_{n}\}$
  	converges to $g$ strongly in $(L^{p^{\prime}}(\Omega))^{N}$ and
  	$(\operatorname{div}g_{n})$ is bounded in $\mathfrak{M}_{b}(\Omega)$; assume also that
  	$\{\rho_{n}\}$ converges to $\mu_{s}^{+}$ and $\{\eta_{n}\}$ to $\mu
  	_{s}^{-}$ in the narrow topology. If $\{u_{n}\}$ is a sequence of renormalized solutions of \eqref{2hvpro1} with data $\mu_n$, 
  	then, up to a subsequence, it converges a.e. in $\Omega$ to a
  	renormalized solution $u$ of problem (\ref{2hvpro1}). Furthermore, $T_{k}(u_{n})$ converges to $T_{k}(u)$  in $W_{0}^{1,p}%
  	(\Omega)$ for any $k>0$.
  \end{theorem}
  
  We also recall the following estimate \cite[Th 2.1]{22PhVe}.
  \begin{proposition}\label{2hvTH4}
  	Let $\Omega$ be a bounded domain of $\mathbb{R}^N$. Then there exists a constant $C>0$, depending on $p$ and $N$ such that if $\mu\in \mathfrak{M}^+_b(\Omega)$ and $u$ is a nonnegative renormalized solution of problem (\ref{2hvpro1}) with  data $\mu$, there holds
  	\begin{equation}
  	\frac{1}{c_{_{86}}}{\bf W}^{\frac{d(x,\partial\Omega)}{3}}_{1,p}[\mu](x) \le u(x)\leq c_{_{86}} {\bf W}^{2\,diam\,(\Omega)}_{1,p}[\mu](x)\qquad\text{a.e. in }\, \Omega.
  	\end{equation}
  \end{proposition} 
  
\nind {\bf Proof of Theorem C.} {\it The condition is necessary}. Assume that \eqref{2hvMT1a} admits a nonnegative renormalized solutions $(u,v)$. By Proposition \ref{2hvTH4} there holds 
  	\begin{align*}
  &	u(x)\geq c_{_{87}} {\bf W}^{\frac{d(x,\partial\Omega)}{3}}_{1,p}[v^{q_1}+\mu](x) \\&
  v(x)\geq c_{_{87}} {\bf W}^{\frac{d(x,\partial\Omega)}{3}}_{1,p}[u^{q_2}+\mu](x) \qquad\text{a.e. in }\, \Omega.
  	\end{align*}
  	Hence,  we infer \eqref{2hvMT1c} from Proposition \ref{260320151}-(ii).\smallskip
  	
\nind {\it Sufficient conditions}.	Let $\{(u_m,v_m)\}_{m\in \mathbb{N}}$ be a 
  	sequence of nonnegative renormalized solutions of the following problems   for $m\in\mathbb{N}$,
  	\begin{equation}
  	\left. \begin{array}{ll}
  	- \Delta _p u_{m+1} = v_m^{q_1}+\mu \qquad&\text{in }\;\Omega  \\ 
  	- \Delta _p v_{m+1} = u_m^{q_2}+\eta \qquad&\text{in }\;\Omega  \\ 
  	\phantom{    - \Delta _p}
  	u_{m+1}=v_{m+1}= 0 \qquad&\text{on }\;\partial\Omega,
  	\end{array} \right.
  	\end{equation}
with initial condition $(u_0,v_0)=0$. The sequences $\{u_m\}$ and $\{v_m\}$ can be constructed in such a way that they are nondecreasing (see e.g. \cite{22PhVe2}).  
  	By Proposition \ref{2hvTH4} we have 
  	\begin{align*}
  &	u_{m+1}\leq c_{_{86}}\mathbf{W}_{1,p}^R[v_m^{q_1}+\mu](x)\\&
  	v_{m+1}\leq c_{_{86}}\mathbf{W}_{1,p}^R[u_m^{q_2}+\eta](x) \qquad\text{a.e. in }\, \Omega,
  	\end{align*}
  	where $R=2\,diam\,(\Omega)$. Thus, by Proposition \ref{pro2} there exists a constant $M_*>0$ depending only on $N,p,q_1,q_2,R$ such that if 
  	\begin{align}
  	\omega(K)\leq M_* \operatorname{Cap}_{\mathbf{G}_{\frac{p(q_1+p-1)}{q_1}},\frac{q_1q_2}{q_1q_2-(p-1)^2}}(K)
  	\end{align}
  	for any compact set $K\subset\mathbb{R}^N$ with $d\omega=\left(\mathbf{W}_{1,p}^{R}[\mu]\right)^{q_2}dx+d\eta$, then 
  	\begin{align}
  	v_m\leq c_{_{73}}  \mathbf{W}_{1,p}^{R}[\omega],~~~u_m\leq c_{_{74}}  \mathbf{W}_{1,p}^{R}[\left(\mathbf{W}^{R}_{1,p}[\omega]\right)^{q_1}]+c_{_{72}} \mathbf{W}_{1,p}^{R}[\mu]~~\forall~~k\geq 0
  	\end{align}
  	in $\Omega$, and 
  	\begin{align}\label{L^q-es}
  	\mathbf{W}_{1,p}^{R}[\omega]\in L^{q_2}(\Omega),~~~~ \mathbf{W}_{1,p}^{R}[\left(\mathbf{W}^{R}_{1,p}[\omega]\right)^{q_1}]+ \mathbf{W}_{1,p}^{R}[\mu]\in L^{q_1}(\Omega).
  	\end{align}
  	This implies that $\{u_m\}, \{v_m\}_{m\in\BBN}$ are well defined and nondecreasing. Thus  $\{(u_m,v_m)\}$ converges a.e in $\Omega$ to some functions $(u,v)$ which satisfies \eqref{thm2-es-uv}  in $\Omega$. Furthermore, we deduce  from \eqref{L^q-es} and the monotone convergence theorem that $u_m^{q_1}\to u^{q_1}$ and $v_m^{q_2}\to u^{q_2}$ in $L^{1}(\Omega)$. Finally we infer that $u$ is a renormalized solution of \eqref{2hvMT1a} by Theorem \ref{2hvP3}.\qeda
	

  \section{p-superharmonic functions and quasilinear equations in $\mathbb{R}^N$ }
  We recall some definitions and properties   of $p$-superharmonic functions (see e.g. \cite {22HeKiMa}, \cite {22KiMa1}, \cite {22KiMa2} for general properties and \cite {VeQ} for a simple presentation). 
  
  \begin{definition} A function $u$ is said to be $p$-harmonic in $\mathbb{R}^N$ if $u\in W^{1,p}_{loc}(\mathbb{R}^N)$ and $-\Delta_p u=0$ in $\mathcal{D}'(\mathbb{R}^N)$; it is always $C^{1}$.
  	A function $u$ is called a $p$-supersolution in $\mathbb{R}^N$ if $u\in W^{1,p}_{loc}(\mathbb{R}^N)$ and $-\Delta_p u\geq 0$ in $\mathcal{D}'(\mathbb{R}^N)$. 
  \end{definition}
  \begin{definition}
  	A lower semicontinuous (l.s.c) function $u: \mathbb{R}^N \to (-\infty,\infty]$ is called $p$-super-harmonic if $u$ is not identically infinite  and if, for all open $D\subset \subset \mathbb{R}^N$ and all $v\in C(\overline{D})$, $p$-harmonic in $D$, $v\le u$ on $\partial D$ implies $v\le u$ in $D$.
  \end{definition}
  Let $u$ be a $p$-superharmonic in $\mathbb{R}^N$. It is well known that $u\wedge k:=\min\{u,k\}\in W^{1,p}_{loc}(\mathbb{R}^N)$ is a p-supersolution for all $k>0$ and $u<\infty$ a.e in $\mathbb{R}^N$, thus, $u$ has a gradient (see the previous section). We also have $|\nabla u|^{p-1}\in L^q_{loc}(\mathbb{R}^N)$, $\frac{|\nabla u|^p}{(|u|+1)^r}\in L^1_{loc}(\mathbb{R}^N)$ and $u\in L^s_{loc}(\mathbb{R}^N)$ for $1\le q<\frac{N}{N-1}$ and $r>1$, $1\le s < \frac{N(p-1)}{N-p}$ (see \cite[Theorem 7.46]{22HeKiMa}).  Thus for any  $0\le \varphi \in C^1_c(\Omega)$, by the dominated convergence theorem, 
  \begin{equation*}
  \left\langle -\Delta_p u,\varphi \right\rangle =\int_{\mathbb{R}^N}|\nabla u|^{p-2}\nabla u.\nabla \varphi dx
  =\mathop {\lim }\limits_{k \to \infty }\int_{\mathbb{R}^N}|\nabla (u\wedge k)|^{p-2}\nabla (u\wedge k) .\nabla \varphi \geq 0.\end{equation*}
  Hence, by the Riesz Representation Theorem, there is a nonnegative Radon measure denoted by $\mu[u]$, called the Riesz measure,  such that $-\Delta_p u=\mu[u]$ in $\mathcal{D}'(\mathbb{R}^N)$.\\
  
  The following weak convergence result for  Riesz measures proved in \cite{22TW4} will be used to obtain the existence of $p$-superharmonic solutions to quasilinear equations. 
  \begin{proposition}\label{2hv07062}
  	Suppose that $\{u_n\}$  is a sequence of nonnegative $p$-superharmonic functions in $\mathbb{R}^N$  that converges a.e to a $p$-superharmonic function $u$. Then the sequence of measures $\{\mu[u_n]\}$ converges to $\mu[u]$ in the weak sense of measures.
  \end{proposition} 
  The proof of the next result can be found in \cite{22PhVe}.
  \begin{proposition}\label{2hv07061}
  	Let $\mu$ be  a measure in $\mathfrak{M}^+(\mathbb{R}^N)$. Suppose that ${\bf W}_{1,p}[\mu]<\infty$ a.e. Then there exists a nonnegative $p$-superharmonic function $u$ in $\mathbb{R}^N$ such that  $-\Delta_p u=\mu $ in $\mathcal{D}'(\mathbb{R}^N)$, $\inf_{\mathbb{R}^N}u=0$ and
  	\begin{equation}\label{2hv07063}
  	\frac{1}{c_{_{86}}}{\bf W}_{1,p}[\mu](x)\le u(x)\le c_{_{86}}{\bf W}_{1,p}[\mu](x),
  	\end{equation}
  	for almost all $x$ in $\mathbb{R}^N$, where the constant $c_{_{86}}$ is the one of Proposition \ref{2hvTH4}.
  	Furthermore any $p$-superharmonic function $u$ in $\mathbb{R}^N$, such that $\inf_{\mathbb{R}^N}u=0$ satisfies \eqref{2hv07063}  with $\mu=-\Delta_pu$. 
  \end{proposition}
  {\bf Proof of Theorem A.} {\it The condition is necessary}.  Assume that \eqref{2hvMT1a} admits a nonnegative $p$-superharmonic functions $(u,v)$. By Proposition \ref{2hv07061} there holds 
  	\begin{align*}
  	&	u(x)\geq c_{_{87}} {\bf W}_{1,p}[v^{q_1}+\mu](x), \\&
  	v(x)\geq c_{_{87}} {\bf W}_{1,p}[u^{q_2}+\eta](x)\qquad \text{for almost all }x\in \Omega.
  	\end{align*}
  	Hence,  we obtain \eqref{2hvMT1cb} from Proposition \ref{260320151}-(i).\medskip
  	
 \nind {\it The condition is sufficient. } Let $\{(u_m,v_m)\}_{m\in \mathbb{N}}$ be a 
  	sequence of nonnegative $p$-superharmonic solutions of the following problems   for $m\in\mathbb{N}$,
  	\begin{equation}
  	\left. \begin{array}{ll}
  	- \Delta _p u_{m+1} = v_m^{q_1}+\mu \qquad&\text{in }\;\mathbb{R}^N  \\ 
  	- \Delta _p v_{m+1} = u_m^{q_2}+\eta \qquad&\text{in }\;\mathbb{R}^N  \\ 
  	\!\!\!\inf_{\mathbb{R}^N}u_{m+1}=\inf_{\mathbb{R}^N}v_{m+1}= 0,
  	\end{array} \right.
  	\end{equation}
  with $(u_0,v_0)=(0,0)$. As in the proof of Theorem C we can assume that $\{u_m\}$ and $\{v_m\}$ are nondecreasing.  
  	By Proposition \ref{2hv07061} we have 
  	\begin{align*}
  	&	u_{m+1}\leq c_{_{86}} \mathbf{W}_{1,p}[v_m^{q_1}+\mu](x)\\&
  	 v_{m+1}\leq c_{_{86}} \mathbf{W}_{1,p}[u_m^{q_2}+\eta](x)\quad\text{for all }\,x\in \Omega.
  	\end{align*}
  Thus, by Proposition \ref{pro1} there exists a constant $c>0$ depending only on $N,p,q_1,q_2$ such that, if 
  	\begin{align}
  	\omega(K)\leq M^* \operatorname{Cap}_{\mathbf{I}_{\frac{p(q_1+p-1)}{q_1}},\frac{q_1q_2}{q_1q_2-(p-1)^2}}(K)
  	\end{align}
  	for any compact set $K\subset\mathbb{R}^N$ with $d\omega=\left(\mathbf{W}_{1,p}[\mu]\right)^{q_2}dx+d\eta$, then there holds in $\Gw$,
  	\begin{align}
  	v_m\leq c_{_{69}} \mathbf{W}_{1,p}[\omega],~~~u_m\leq c_{_{70}} \mathbf{W}_{1,p}[\left(\mathbf{W}_{1,p}[\omega]\right)^{q_1}]+c_{_{68}} \mathbf{W}_{1,p}[\mu]\quad\text{for all }\,m\geq 0,
  	\end{align}
  	 and 
  	\begin{align}\label{L^q-es'}
  	\mathbf{W}_{1,p}[\omega]\in L^{q_2}_{loc}(\mathbb{R}^N),~~~~ \mathbf{W}_{1,p}[\left(\mathbf{W}_{1,p}[\omega]\right)^{q_1}]+ \mathbf{W}_{1,p}[\mu]\in L^{q_1}_{loc}(\mathbb{R}^N).
  	\end{align}
  	This implies that $\{u_m\}, \{v_m\}$ are well defined and nondecreasing. Thus  $\{(u_m,v_m)\}$ converges a.e in $\mathbb{R}^N$ to some functions $(u,v)$ which satisfies \eqref{thm2-es-uv}  in $\mathbb{R}^N$. Furthermore, we infer  from \eqref{L^q-es} and the monotone convergence theorem that $u_m^{q_1}\to u^{q_1}, v_m^{q_2}\to u^{q_2}$ in $L^{1}_{loc}(\mathbb{R}^N)$. By Proposition \ref{2hv07062} we deduce that $(u,v)$ are nonnegative $p$-superharmonic solutions of \eqref{2hvMT1ab}.\phantom{--------}
  \qeda
  

   \section{Hessian equations}
   In this section $\Omega\subset\mathbb{R}^N$ is either a bounded domain with a $C^2$ boundary or the whole $\mathbb{R}^N$. For $k=1,...,N$ and $u\in C^2(\Omega)$ the k-hessian operator $F_k$ is defined by 
   $$F_k[u]=S_k(\lambda(D^2u)),$$
   where $\lambda(D^2u)=\lambda=(\lambda_1,\lambda_2,...,\lambda_N)$ denotes the eigenvalues of the Hessian matrix of second partial derivative $D^2u$ and $S_k$ is the k-th elementary symmetric polynomial that is \[{S_k}(\lambda ) = \sum\limits_{1 \le {i_1} < ... < {i_k} \le N} {{\lambda _{{i_1}}}...{\lambda _{{i_k}}}}. \]
   We can see that \[{F_k}[u] = {\left[ {{D^2}u} \right]_k},\]
   where for a matrix $A=(a_{ij})$, $[A]_k$ denotes the sum of the k-th principal minors.
   We assume that  $\partial \Omega$ is uniformly (k-1)-convex, that is $$S_{k-1}(\kappa ) \geq c_0 >0~on ~~ \partial\Omega,$$
   for some positive constant $c_0$, where $\kappa= (\kappa_1,\kappa_2,...,\kappa_{n-1})$ denote the principal curvatures of $\partial \Omega$ with respect to its inner normal.
   \begin{definition}\label{2hvk-conv}
   	An upper-semicontinuous function $u:\Omega\to [-\infty,\infty)$ is k-convex (k-subharmonic) if, for every open set $\Omega'\subset\overline \Omega'\subset\Omega$ and for every function $v\in C^2(\Omega')\cap C(\overline{\Omega'})$ satisfying $F_k[v]\leq 0$ in $\Omega'$, the following implication is true 
   	$$u\leq v ~on ~\partial\Omega' ~~~\Longrightarrow~~~ u\leq v~~in ~~\Omega'.$$
   	We denote by $\Phi^k(\Omega)$ the class of all $k$-subharmonic functions in $\Omega$ which are not identically equal to $-\infty$.   
   \end{definition} 
   The following weak convergence result for $k$-Hessian operators proved in \cite{22TW2} is fundamental in our study. 
   \begin{proposition}\label{2hvTH6} Let $\Omega$ be either a bounded uniformly (k-1)-convex in $\mathbb{R}^N$ or the whole  $\mathbb{R}^N$. For each $u\in \Phi^k(\Omega)$, there exists a nonnegative Radon measure $\mu_k[u]$ in $\Omega$ such that \smallskip
   	
   	\noindent {\bf 1} $\mu_k[u]=F_k[u]$ for $u\in C^2(\Omega)$.\smallskip
   	
   	\noindent {\bf 2}  If $\{u_n\}$ is a sequence of k-convex functions which converges a.e to $u$, then $\mu_{k}[u_n]\rightharpoonup\mu_{k}[u] $ in the weak sense of measures.
   	
   \end{proposition}
   
   As in the case of quasilinear equations with measure data, precise estimates of solutions of k-Hessian equations with measures data are expressed in terms of Wolff potentials. The next results are proved in \cite{22TW2,22La,22PhVe}. 
   \begin{theorem}\label{2hvTH5}Let $\Omega\subset \mathbb{R}^N$ be a bounded $C^2$, uniformly (k-1)-convex  domain. Let  $\mu$ be a nonnegative Radon measure in $\Gw$ which can be decomposed under the form 
   	$$\mu=\mu_1+f,$$
   	where $\mu_1$ is a measure with compact support in $\Omega$ and $f\in L^q(\Omega)$ for some $q>\frac{N}{2k}$ if $k\le\frac{N}{2}$, or $p=1$ if $k>\frac{N}{2}$. Then there exists a nonnegative function $u$ in $\Omega$, continuous near $\partial \Omega$, such that $-u\in \Phi^k(\Omega)$  and u is a solution of the problem 
   	\[\begin{array}{lll}
   	{F_k}[-u] = \mu \qquad&\;in\;\Omega,  \\ 
   	\phantom{{F_k}[-]}u =0\qquad&\;on\;\partial\Omega.
   	\end{array} \]
   	Furthermore,  any nonnegative function $u$ such that  $-u\in \Phi^k(\Omega)$ which is continuous near $\partial \Omega$ and is a solution of above equation, satisfies 
   	\begin{equation}
   	\frac{1}{c_{_{88}}}{\bf W}^{\frac{d(x,\partial \Omega)}{8}}_{\frac{2k}{k+1},k+1}[\mu]\le u(x)\le c_{_{88}}  {\bf W}^{2diam\,\Omega}_{\frac{2k}{k+1},k+1}[\mu](x),
   	\end{equation}
   	where  $c_{_{88}}$ is a positive constant independent of $x,u$ and $\Omega$. 
   \end{theorem}
   \begin{theorem}\label{2hv07065}
   	Let $\mu$ be  a measure in $\mathfrak{M}^+(\mathbb{R}^N)$ and $2k<N$. Suppose that ${\bf W}_{\frac{2k}{k+1},k+1}[\mu]<\infty$ a.e. Then there exists $u$, $-u\in \Phi^k(\mathbb{R}^N)$ such that $\inf_{\mathbb{R}^N}u=0$ and $F_k[-u]=\mu ~~\text{ in }~~ \mathbb{R}^N$
   	and 
   	\begin{equation}\label{2hv07064}
   	\frac{1}{c_{_{88}}}{\bf W}_{\frac{2k}{k+1},k+1}[\mu](x)\le u(x)\le c_{_{88}}{\bf W}_{\frac{2k}{k+1},k+1}[\mu](x),
   	\end{equation}
   	for all $x$ in $\mathbb{R}^N$.
   	Furthermore, if $u$ is a nonnegative function  such that $\inf_{\mathbb{R}^N}u=0$ and  $-u\in \Phi^k(\mathbb{R}^N)$, then \eqref{2hv07064} holds  with $\mu=F_k[-u]$. 
   \end{theorem}
   

\nind{\bf Proof of Theorem E.} {\it The condition is necessary}. Assume that \eqref{2hvMT3a}  admits a nonnegative solution $(u,v)$, continuous near $\partial \Omega$, such that $-u,-v\in \Phi^k(\Omega)$ and $u^{s_2},v^{s_1}\in L^1(\Omega)$. Then by Theorem \ref{2hvTH5} we have 
	\begin{align*}
&	u(x)\geq \frac{1}{c_{_{88}}}{\bf W}^{\frac{d(x,\partial\Omega)}{8}}_{\frac{2k}{k+1},k+1}[v^{s_1}+\mu](x)\\&
	v(x)\geq \frac{1}{c_{_{88}}}{\bf W}^{\frac{d(x,\partial\Omega)}{8}}_{\frac{2k}{k+1},k+1}[u^{s_2}+\eta](x)
\qquad\textrm{ for almost all } x\in \Omega.
	\end{align*}
	Using the part 2 of Proposition \ref{260320151}, we conclude that \eqref{2hvMT3c} holds.\smallskip
	
\nind{\it The condition is sufficient}. We define  a sequence of nonnegative functions $u_m,v_m$, continuous near $\partial \Omega$ and such that $-u_m,-v_m\in \Phi^k(\Omega)$, by the following iterative scheme  for $m\geq 0$,
	\begin{equation}\label{2hvk-hess-1}
	\begin{array}{ll}
	F_k[-u_{m+1}]  =v_m^{s_1}+ \mu\qquad&\text{in }\;\Omega,   \\ 
		F_k[-v_{m+1}]  =u_m^{s_2}+ \eta\qquad&\text{in }\;\Omega,   \\ 
	\phantom{ F_k[-]}
	u_{m+1} =v_{m+1}=0&\text{on }\;\partial\Omega.
	\end{array}
	\end{equation}
	Clearly, we can assume that $\{u_m\}$ is nondecreasing as in \cite{22PhVe2}. By Theorem \ref{2hvTH5} we have 
	\begin{equation}\label{2hvk-hess-2}  \begin{array}{ll}
	u_{m+1}\leq c_{_{88}} {\bf W}^{R}_{\frac{2k}{k+1},k+1}[v_m^{s_1}+
	\mu]\;,\quad
v_{m+1}\leq c_{_{88}} {\bf W}^{R}_{\frac{2k}{k+1},k+1}[u_m^{s_2}+
	\mu]\qquad\text{in}~\Omega,
	\end{array}   \end{equation}
	where $R=2 \,diam\, (\Omega)$.\\
Then, by Proposition \ref{pro2},  there exists a constant $M_*>0$ depending only on $N,p,q_1,q_2,R$ such that if 
	\begin{align*}
	\omega(K)\leq M_* \operatorname{Cap}_{\mathbf{G}_{\frac{2k(s_1+k)}{s_1},\frac{s_1s_2}{s_1s_2-k^2}}}(K)
	\end{align*}
	for any compact set $K\subset\mathbb{R}^N$ with $d\omega=\left({\bf W}^{R}_{\frac{2k}{k+1},k+1}\mu]\right)^{s_2}dx+d\eta$, then there holds, 
	\begin{align*}
	v_m\leq c_{_{73}} {\bf W}^{R}_{\frac{2k}{k+1},k+1}[\omega],~~~u_m\leq c_{_{74}} {\bf W}^{R}_{\frac{2k}{k+1},k+1}[\left({\bf W}^{R}_{\frac{2k}{k+1},k+1}[\omega]\right)^{s_1}]+c_{_{72}} {\bf W}^{R}_{\frac{2k}{k+1},k+1}[\mu]
	\end{align*}
	in $\Omega$, for all $m\in\BBN$, for some positive constants $c_{_{72}},c_{_{73}}$ and $c_{_{74}}$ depending only on $N,k,s_1,s_2,R$. 
	Note that we can write 
	$$v_m^{s_1}+\mu = \left(\mu_1 + \chi_{_{\Omega_\delta}}v_m^{s_1}\right)+\left((1-\chi_{_{\Omega_\delta}})v_m^{s_1}+f\right), $$  
	and 
	$$u_m^{s_2}+\eta = \left(\eta_1 + \chi_{_{\Omega_\delta}}u_m^{s_2}\right)+\left((1-\chi_{_{\Omega_\delta}})u_m^{s_2}+g\right), $$  
	where $\Omega_\delta=\{x\in\Omega:d(x,\partial \Omega)>\delta\}$ and  $\delta>0$ is small enough and  since $u_m$ is continuous near $\partial \Omega$, then  $v_m^{s_1}+\mu, u_m^{s_2}+\eta $ satisfy the assumptions of the data  in Theorem \ref{2hvTH5}.   Therefore the sequence $\{u_m\}$ is well defined and nondecreasing. Thus, $\{u_m\}$ converges a.e in $\Omega$ to some function $u$ which satisfies \eqref{thm3-eswolff} in $\Omega$. Furthermore, by the monotone convergence theorem there holds $v_m^{s_1}\to v,u_m^{s_2}\to u$ in $L^1(\Omega)$. Finally, by Proposition \ref{2hvTH6}, we infer  that \eqref{2hvMT3a} admits a nonnegative  solutions $u,v$, continuous near $\partial \Omega$, with $-u,-v\in \Phi^k(\Omega)$ satisfying \eqref{thm3-eswolff}.
\qeda
\medskip

\nind{\bf Proof of Theorem F} {\it The condition is necessary}.  Assume that \eqref{2hvMT3a}  admits nonnegative solution $(u,v)$, such that $-u,-v\in \Phi^k(\mathbb{R}^N)$ and $u^{s_2},v^{s_1}\in L^1_{loc}(\mathbb{R}^N)$. Then by Theorem \ref{2hvTH5} we have 
	\begin{align*}
	&	u(x)\geq \frac{1}{c_{_{88}}}{\bf W}_{\frac{2k}{k+1},k+1}[v^{s_1}+\mu](x)\\&
	v(x)\geq \frac{1}{c_{_{88}}}{\bf W}_{\frac{2k}{k+1},k+1}[u^{s_2}+\eta](x)
	\qquad\textrm{ for almost all } x\in \mathbb{R}^N.
	\end{align*}
	Using  Proposition \ref{260320151}-(ii), we conclude that \eqref{2hvMT3c} holds.\smallskip
	
\nind {\it The condition is sufficient}.  We defined  a sequence of nonnegative functions $u_m,v_m$, continuous near $\partial \Omega$ and such that $-u_m,-v_m\in \Phi^k(\Omega)$, by the following iterative scheme  for $m\geq 0$,
	\begin{equation*}
	\begin{array}{ll}
	\!F_k[-u_{m+1}]  =v_m^{s_1}+ \mu\qquad&\text{in }\;\mathbb{R}^N,   \\ 
	F_k[-v_{m+1}]  =u_m^{s_2}+ \eta\qquad&\text{in }\;\mathbb{R}^N,   \\ 
	\phantom{}
	\!\!\inf_{\mathbb{R}^N}u_{m+1} =\inf_{\mathbb{R}^N}v_{m+1}=0.
	\end{array}
	\end{equation*}
	As in the previous proofs $\{u_m\}$ is nondecreasing. By Theorem \ref{2hvTH5} we have 
	\begin{equation*} \begin{array}{ll}
	u_{m+1}\leq c_{_{88}} {\bf W}_{\frac{2k}{k+1},k+1}[v_m^{s_1}+
	\mu]\\
	v_{m+1}\leq c_{_{88}} {\bf W}_{\frac{2k}{k+1},k+1}[u_m^{s_2}+
	\mu]\qquad\text{a.e. in}~\mathbb{R}^N.
	\end{array}   \end{equation*}
	Then, by Proposition \ref{pro1},  there exists a constant $M^*>0$ depending only on $N,p,q_1,q_2,R$ such that if 
	\begin{align*}
	\omega(K)\leq M^* \operatorname{Cap}_{\mathbf{I}_{\frac{2k(s_1+k)}{s_1},\frac{s_1s_2}{s_1s_2-k^2}}}(K)
	\end{align*}
	for any compact set $K\subset\mathbb{R}^N$ with $d\omega=\left({\bf W}_{\frac{2k}{k+1},k+1}\mu]\right)^{s_2}dx+d\eta$, then 
	\begin{align*}
	v_m\leq c_{_{69}} {\bf W}_{\frac{2k}{k+1},k+1}[\omega],~~~u_m\leq c_{_{70}} {\bf W}_{\frac{2k}{k+1},k+1}[\left({\bf W}_{\frac{2k}{k+1},k+1}[\omega]\right)^{s_1}]+c_{_{68}} {\bf W}_{\frac{2k}{k+1},k+1}[\mu]
	\end{align*}
	in $\Omega$, for all $m\in\BBN$, where $c_{_{68}}$, $c_{_{69}}$ and $c_{_{70}}$ depend on $N,k,s_1,s_2,R$. 
 Therefore the sequence $\{u_m\}$ is well defined and nondecreasing. Thus, $\{u_m\}$ converges a.e in $\Omega$ to some function $u$ for which \eqref{thm4-eswolff} is satisfied in $\mathbb{R}^N$. Furthermore, by the monotone convergence theorem we have $v_m^{s_1}\to v,u_m^{s_2}\to u$ in $L^1_{loc}(\mathbb{R}^N)$. Finally, by Proposition \ref{2hvTH6}, we obtain  that \eqref{2hvMT3a} 	admits a nonnegative  solutions $u,v$ with $-u,-v\in \Phi^k(\mathbb{R}^N)$ satisfying 
	\eqref{thm4-eswolff}.
\qeda
\section{Further results}
The method exposed in the previous sections, can be applied to types of problems. We give below an example for a semilinear system in  
$\BBR^N_+=\{x=(x',x_{_N}), x'\in \mathbb{R}^{N-1},x_{_N}>0\}$.
\begin{equation}\label{0411201410d*}
\begin{array}
{llll}
-\Delta u=v^{q_1}\qquad&\text{in} ~\mathbb{R}^N_+\\[1mm]
-\Delta v=u^{q_2}\qquad&\text{in} ~\mathbb{R}^N_+\\[1mm]  
\phantom{-\Delta}
u=\sigma_1,v=\sigma_2\qquad&\text{on} ~\prt\mathbb{R}^N_+\approx \mathbb{R}^{N-1},                                                                               
\end{array}
\end{equation}
where we have identified $\prt\BBR^N_+$ and $\BBR^{N-1}$. We denote by $\mathbf{P}$ (resp. $\mathbf{G}$) the Poisson kernel  in $\mathbb{R}^N_+$ (resp the Green kernel in $\mathbb{R}^N$). The Poisson potential and the Green potential, $\mathbf{P}[.]$ and $\mathbf{G}[.]$, associated to $-\Gd$   are defined respectively by 
\begin{align*}
\mathbf{P}[\sigma](y)=\int_{\partial\mathbb{R}_+^N}{P}(y,z)d\sigma(z), ~~{G}[f](y)=\int_{\mathbb{R}_+^N}\mathbf{G}(y,x)f(x)dx,
\end{align*}
see \cite{MV5}. We set $\rho(x)=x_{_N}$ and define the capacity $\operatorname{Cap}_{\alpha,s}^\rho$ by 
\begin{align*}
\operatorname{Cap}_{\alpha,s}^\rho(K)=\inf\left\{\int_{\mathbb{R}^N_+} f^s\rho\,dx: f\geq 0, \mathbf{I}_{\alpha}[f\rho\,\chi_{_{\mathbb{R}^N_+}}]\geq \chi_{_K}\right\},
\end{align*}
for all Borel set $K\subset\BBR^N$, where $\mathbf{I}_{\alpha}$ is the Riesz kernel of order $\ga$ in $\BBR^N$.
\begin{theorem} \label{traceproblem}Let $1\leq q_1<\frac{N}{N-1}$, $q_1q_2>1$.  If there exists a constant $\tilde c>0$ such that if 
\bel{sys-1}\BA {lll}
(i)\qquad\qquad	&\myint{K}{}\rho(x) (\mathbf{P}[\sigma_1](x))^{q_2}dx\leq \tilde c \operatorname{Cap}_{\frac{q_1+2}{q_1},\frac{q_1q_2}{q_1q_2-1}}^\rho(K),	\qquad\qquad\qquad\qquad\\[4mm]
(ii)\qquad	&\sigma_2(G)\leq \tilde c \operatorname{Cap}_{I_{\frac{2(q_2+1)}{q_1q_2},\frac{q_1q_2}{q_1q_2-1}}}(G),
\EA	\ee
for all Borel sets $K\subset\BBR_+^N$ and $G\subset\BBR^{N-1}$, then the problem \eqref{0411201410d*} admits a solution. 
\end{theorem}
All  solutions in above theorem  are understood in the usual very weak sense:  $u\in L^1_{loc}(\mathbb{R}^N_+\cap B)$, $u^{q_2},v^{q_1}\in L^1_{\gr}(\mathbb{R}^N_+\cap B)$  for any ball $B$  and 
\begin{align*}
&\int_{\mathbb{R}^N_+} u (-\Delta \xi)dx =\int_{\mathbb{R}^N_+} v^{q_1}\xi dx-\int_{\partial\mathbb{R}^N_+}\frac{\partial \xi}{\partial n}d\sigma_1,\\&
\int_{\mathbb{R}^N_+} v (-\Delta \xi)dx =\int_{\mathbb{R}^N_+} u^{q_2}\xi dx-\int_{\partial\mathbb{R}^N_+}\frac{\partial \xi}{\partial n}d\sigma_2,
\end{align*}
for any $\xi \in C^2(\overline{\mathbb{R}^N_+})\cap C_c(\mathbb{R}^N)$ with $\xi=0$ on $\partial\mathbb{R}^N_+$.
It is well-known that  such a solution  $u$ satisfies 
\begin{align*}
u=\mathbf{G}[v^{q_1}]+\mathbf{P}[\sigma_1],~~v=\mathbf{G}[u^{q_2}]+\mathbf{P}[\sigma_2]\qquad\text{a.e. in } \mathbb{R}_+^N.
\end{align*}

To prove Theorem \ref{traceproblem} we need the following basic estimate,
\begin{lemma} \label{lela}Assume that $0< q_1<\frac{N}{N-1}$. Then for any $\gw\in\mathfrak M^+_b(\mathbb{R}^N)$, 
	\begin{align}\label{es1'}
	\mathbf{I}_2\left[\left(\mathbf{I}_1[\omega]\right)^{q_1}\right]\leq 
	c_{_{89}} \mathbf{W}_{\frac{q_1+2}{q_1+1},\frac{q_1+1}{q_1}}[\omega]\qquad\text{a.e. in }\mathbb{R}^N,
	\end{align}
	where $c_{_{89}}>0$ depends on $q_1,q_2$ and $N$.
\end{lemma}
\proof 
The proof of Lemma \ref{lela} is similar to the one of Lemma \ref{241220142} and details are omitted. Note that if $\gw\in\mathfrak M_b(\overline{\mathbb{R}^N_+})$ it is extended by $0$ in $\mathbb{R}^N_-$.
\qeda
\begin{remark}The condition $0<q_1<\frac{N}{N-1}$ is a necessary and sufficient condition in order  $\left(\mathbf{I}_1[\omega]\right)^{q_1}$ be locally integrable in $\mathbb {R}^N$  for any $\gw\in\mathfrak M^+_b(\mathbb{R}^N)$. 
\end{remark}
\begin{theorem} \label{final}Let $ q_1\geq 1$, $q_1q_2>1$ and $\gw\in\mathfrak M_b(\overline{\mathbb{R}^N_+})$. If 	
\begin{align*}
	\omega(K)\leq c_{_{90}} \operatorname{Cap}_{\frac{q_1+2}{q_1},\frac{q_1q_2}{q_1q_2-1}}^\rho(K)\qquad \forall~ K\subset \overline{\mathbb{R}^N_+},\;
	K\text{ Borel},
	\end{align*}
for some $c_{_{90}}>0$,	then
	\begin{align}\label{2412201410}
	\mathbf{I}_1\left[\left(\mathbf{W}_{\frac{q_1+2}{q_1+1},\frac{q_1+1}{q_1}}[\omega]\right)^{q_2}\rho\,\chi_{_{\mathbb{R}^N_+}}\right]\leq c_{_{91}} \mathbf{I}_1[\omega]\qquad\text{a.e. in }\mathbb{R}^N_+.
	\end{align}
\end{theorem}
\proof {\it Step 1}.  For any compact $K\subset \left\{x\in\BBR^N_+:\mathbf{I}_{\frac{q_1+2}{q_1}}[f\rho\,\chi_{_{\mathbb{R}^N_+}}](x)>\lambda\right\}$, we have 
	\begin{align*}
	\omega(K)\leq c_{_{90}} \operatorname{Cap}_{\frac{q_1+2}{q_1},\frac{q_1q_2}{q_1q_2-1}}^\rho(K)
	\leq c_{_{90}} \lambda^{-\frac{q_1q_2}{q_1q_2-1}}\int_{\mathbb{R}^N_+} f^{\frac{q_1q_2}{q_1q_2-1}}\rho\, dx
	\end{align*}
	by assumption and the definition of the capacity. Hence, 
	\begin{align*}
	\lambda^{\frac{q_1q_2}{q_1q_2-1}} \omega\left(\left\{\mathbf{I}_{\frac{q_1+2}{q_1}}[f\rho\,\chi_{_{\mathbb{R}^N_+}}]>\lambda\right\}\right)\leq c_{_{90}} \int_{\mathbb{R}^N_+} f^{\frac{q_1q_2}{q_1q_2-1}}\rho\, dx\qquad\forall~\lambda>0.
	\end{align*}
	This implies an estimate in Lorentz space, 
	\begin{align}\label{241220148}
	||\mathbf{I}_{\frac{q_1+2}{q_1}}[f\rho\,\chi_{_{\mathbb{R}^N_+}}]||_{L^{\frac{q_1q_2}{q_1q_2-1},\infty}(\mathbb{R}^N,d\omega)}\leq ||f||_{L^{\frac{q_1q_2}{q_1q_2-1}}(\mathbf{R}^N,\chi_{_{\mathbb{R}^N_+}}\rho\, dx)} \qquad\forall~f\geq 0.
	\end{align}\smallskip
	
\nind	{\it Step 2}.  Since, for any $g\in C_c(\mathbb{R}^N_+)$, 
	\begin{align*}
	\int_{\mathbb{R}^N_+}\mathbf{I}_{\frac{q_1+2}{q_1}}[g\omega]f\rho\, dx=\int_{\mathbb{R}^N}\mathbf{I}_{\frac{q_1+2}{q_1}}[f\rho\, \chi_{_{\mathbb{R}^N_+}}]gd\omega,
	\end{align*}
we infer, using duality  between $L^{p,1}$ and $L^{p',\infty}$, Holder's inequality therein and \eqref{241220148}, that 
	\begin{align*}
	\int_{\mathbb{R}^N_+}\mathbf{I}_{\frac{q_1+2}{q_1}}[g\omega]f\rho\, dx&\leq ||\mathbf{I}_{\frac{q_1+2}{q_1}}[f\rho\,\chi_{_{\mathbb{R}^N_+}}]||_{L^{\frac{q_1q_2}{q_1q_2-1},\infty}(\mathbb{R}^N,d\omega)}||g||_{L^{q_1q_2,1}(\mathbb{R}^N,d\omega)}\\& \leq 
	 ||f||_{L^{\frac{q_1q_2}{q_1q_2-1}}(\mathbb{R}^N,\chi_{_{\mathbb{R}^N_+}}\rho\, dx)} ||g||_{L^{q_1q_2,1}(\mathbb{R}^N,d\omega)} \qquad\forall~f,g\geq 0.
	\end{align*}
	Therefore, 
	\bel{E1}
	||\mathbf{I}_{\frac{q_1+2}{q_1}}[g\omega]||_{L^{q_1q_2}(\mathbb{R}^N,\chi_{_{\mathbb{R}^N_+}}\rho\, dx)}\leq ||g||_{L^{q_1q_2,1}(\mathbb{R}^N,d\omega)}.  
	\ee\smallskip
	
\nind	{\it Step 3}.    Taking $g=\chi_{_{B_t(x)}}$ and since for $q_1\geq 1$
\begin{align*}
\mathbf{W}_{\frac{q_1+2}{q_1+1},\frac{q_1+1}{q_1}}[\nu](x)&=\int_{0}^{\infty}\left(\frac{\nu(B_\rho(x))}{\rho^{N-\frac{q_1+2}{q_1}}}\right)^{q_1}dx\\&~~\leq c_{_{89}} \left(\int_{0}^{\infty}\frac{\nu(B_\rho(x))}{\rho^{N-\frac{q_1+2}{q_1}}}dx\right)^{q_1}\\&= c_{_{89}}
\left(\mathbf{I}_{\frac{q_1+2}{q_1}}[\nu](x)\right)^{q_1}\qquad\forall\gn\in\mathfrak M_b^+(\BBR^N),~~\forall x\in \mathbb{R}^N,
\end{align*}
 we deduce that for almost all $x\in \BBR^N_+$,
	\begin{align*}
	\int_{\mathbb{R}^N_+}\left(\mathbf{W}_{\frac{q_1+2}{q_1+1},\frac{q_1+1}{q_1}}[\chi_{_{B_t(x)}}\omega]\right)^{q_2}\rho\, dy\leq c_{_{90}}\omega(B_{t}(x)),
	\end{align*}
from (\ref{E1}), which implies
\begin{align}\label{ab}
	\omega(B_t(x))\leq c_{_{91}} \frac{t^{(N-\frac{q_1+2}{q_1})\frac{q_1q_2}{q_1q_2-1}}}{\left(\myint{B_{2t}(x)}{}\chi_{_{\mathbb{R}^N_+}}\rho\, dy\right)^{\frac{1}{q_1q_2-1}}}\leq c_{_{92}} \frac{t^{n-\frac{q_1+2}{q_1}\frac{q_1q_2}{q_1q_2-1}}}{(\max\{x_n,t\})^{\frac{1}{q_1q_2-1}}},
	\end{align}
	since $\int_{B_r(x)}\chi_{_{\mathbb{R}^N_+}}\rho\, dy \asymp r^N\max\{x_{_N},r\}$ for any $x\in \mathbb{R}^N_+,r>0$ where the symbol $\asymp$ is defined by 
	$$A\asymp B\Longleftrightarrow \myfrac {1}{c}B\leq A\leq cB\quad\text{for some constant }c>0.$$
It implies also 
	\begin{align}\label{241220149}
	\myint{B_t(x)}{}\left(\mathbf{W}^t_{\frac{q_1+2}{q_1+1},\frac{q_1+1}{q_1}}[\omega]\right)^{q_2}\chi_{_{\mathbb{R}^N_+}}\rho\, dy\leq c_{_{90}}\omega(B_{2t}(x)),
	\end{align}
from which follows  
	\begin{align*}
	\int_{0}^{\infty}\frac{1}{t^N}\myint{B_t(x)}{}\left(\mathbf{W}^t_{\frac{q_1+2}{q_1+1},\frac{q_1+1}{q_1}}[\omega]\right)^{q_2}\chi_{_{\mathbb{R}^N_+}}\rho\, dydt\leq c_{_{93}} \mathbf{I}_1[\omega](x).
	\end{align*}
	Therefore, if  the following inequality holds 
	\begin{align}\label{2412201411}
	\int_{0}^{\infty}\frac{1}{t^N}\myint{B_t(x)}{}\left(\int_{t}^{\infty}\left(\frac{\omega(B_r(y))}{r^{N-\frac{q_1+2}{q_1}}}\right)^{q_1}\frac{dr}{r}\right)^{q_2}\chi_{_{\mathbb{R}^N_+}}\rho\, dydt\leq c_{_{93}} \mathbf{I}_1[\omega](x),
	\end{align}
it will imply \eqref{2412201410}. \smallskip
	
\nind	{\it Step 4}.  We claim that \eqref{2412201411} holds.
	Since $B_r(y)\in B_{2r}(x)$, $y\in B_t(x),r\geq t$, 
	\begin{align*}
	&\int_{0}^{\infty}\frac{1}{t^N}\myint{B_t(x)}{}\left(\int_{t}^{\infty}\left(\frac{\omega(B_r(y))}{r^{N-\frac{q_1+2}{q_1}}}\right)^{q_1}\frac{dr}{r}\right)^{q_2}\chi_{_{\mathbb{R}^N_+}}\rho\, dydt\\&\phantom{-----}\leq 
	\int_{0}^{\infty}\frac{1}{t^N}\myint{B_t(x)}{}\chi_{_{\mathbb{R}^N_+}}\rho\,dy\left(\int_{t}^{\infty}\left(\frac{\omega(B_{2r}(x))}{r^{N-\frac{q_1+2}{q_1}}}\right)^{q_1}\frac{dr}{r}\right)^{q_2}dt\\&\phantom{-----}\asymp  \int_{0}^{\infty}\max\{x_n,t\}\left(\int_{t}^{\infty}\left(\frac{\omega(B_{2r}(x))}{r^{N-\frac{q_1+2}{q_1}}}\right)^{q_1}\frac{dr}{r}\right)^{q_2}dt.
	\end{align*}
	By integration by part,
	\begin{align*}
	&\int_{0}^{\infty}\frac{1}{t^N}\myint{B_t(x)}{}\left(\int_{t}^{\infty}\left(\frac{\omega(B_r(y))}{r^{N-\frac{q_1+2}{q_1}}}\right)^{q_1}\frac{dr}{r}\right)^{q_2}\chi_{_{\mathbb{R}^N_+}}\rho\,dydt \\& =q_2\int_{0}^{\infty}\!\!\!\int_{0}^{t}\max\{x_{_N},s\}ds\left(\int_{t}^{\infty}\left(\frac{\omega(B_{2r}(x))}{r^{N-\frac{q_1+2}{q_1}}}\right)^{q_1}\frac{dr}{r}\right)^{q_2-1}\!\!\!\left(\frac{\omega(B_{2t}(x))}{t^{N-\frac{q_1+2}{q_1}}}\right)^{q_1}\frac{dt}{t}
	\\& =q_2\int_{0}^{\infty}\!\!\!\int_{0}^{t}\max\{x_{_N},s\}ds\left(\int_{t}^{\infty}\left(\frac{\omega(B_{2r}(x))}{r^{N-\frac{q_1+2}{q_1}}}\right)^{q_1}\frac{dr}{r}\right)^{q_2-1}\!\!\!\left(\frac{\omega(B_{2t}(x))}{t^{N-\frac{q_1+2}{q_1}}}\right)^{q_1-1}\!\!\!t^{\frac2{q_1}}\frac{\omega(B_{2t}(x))}{t^{N-1}}\frac{dt}{t}
	\end{align*}
	We have 
	\begin{align*}
	\int_{0}^{t}\max\{x_{_N},s\}ds\asymp t\max\{x_{_N},t\},
	\end{align*}
	\begin{align*}
	\left(\int_{t}^{\infty}\left(\frac{\omega(B_{2r}(x))}{r^{N-\frac{q_1+2}{q_1}}}\right)^{q_1}\frac{dr}{r}\right)^{q_2-1}&\leq c_{_{94}}\left(\int_{t}^{\infty}\left(\frac{r^{-\frac{q_1+2}{q_1(q_1q_2-1)}}}{(\max\{x_{_N},r\})^{\frac{1}{q_1q_2-1}}}\right)^{q_1}\frac{dr}{r}\right)^{q-1}\\& \asymp t^{-\frac{(q_1+2)(q_2-1)}{q_1q_2-1}}(\max\{x_{_N},t\})^{-\frac{q_1(q_1-1)}{q_1q_2-1}},
	\end{align*}
	by (\ref{ab}) and 
	\begin{align*}
	\left(\frac{\omega(B_{2t}(x))}{t^{N-\frac{q_1+2}{q_1}}}\right)^{q_1-1}t^{\frac{2}{q_1}}&\leq c_{_{95}}\left(\frac{t^{-\frac{q_1+2}{q_1(q_1q_2-1)}}}{(\max\{x_{_N},t\})^{\frac{1}{q_1q_2-1}}}\right)^{q_2-1}t^{\frac{2}{q_1}}\\&=c_{_{95}} t^{-\frac{(q_1+2)(q_1-1)}{q_1(q_1q_2-1)}+\frac{2}{q_1}}(\max\{x_{_N},t\})^{-\frac{q_1-1}{q_1q_2-1}}.
	\end{align*}
	Thus, 
	\begin{align*}
	\int_{0}^{t}\max\{x_{_N},s\}ds\left(\int_{t}^{\infty}\left(\frac{\omega(B_{2r}(x))}{r^{N-\frac{q_1+2}{q_1}}}\right)^{q_1}\frac{dr}{r}\right)^{q_2-1}\left(\frac{\omega(B_{2t}(x))}{t^{N-\frac{q_1+2}{q_1}}}\right)^{q_1-1}t^{2/q_1}\leq c_{_{96}},
	\end{align*}
	and we obtain \eqref{2412201411}. 
\qeda

\begin{lemma} Let $\ga>0$, $s>1$ such that $\ga+\frac{2}{s'}<N-1$ where $s'=\frac s{s-1}$. For all $\eta\in \mathfrak{M}^+(\mathbb{R}^{N-1})$, there holds 
	\begin{align}\label{es8}
	\int_{\mathbb{R}^N}(\mathbf{I}_\alpha[\eta\otimes \delta_{\{x_{_N}=0\}}] )^{s'}x_{_N}dx\asymp \int_{\mathbb{R}^{N-1}}\left(\int_{0}^{\infty}\frac{\eta(B'_t(x'))}{t^{N-1-\alpha-\frac{2}{s'}}}\frac{dt}{t}\right)^{s'}dx',
	\end{align}
where $I_\beta$ is the Riesz potential of order $\beta$ in $\mathbb{R}^{N-1}$. As a consequence, we have 
	\begin{align}\label{es9}
	\operatorname{Cap}_{\alpha,s}^\rho(E\times\{x_{_N}=0\})\asymp \operatorname{Cap}_{I_{\alpha+2/s'-1},s}(E)\qquad\forall E\subset \mathbb{R}^{N-1}\,,\; E\text { Borel}.
	\end{align}
\end{lemma}
\proof We have 
	\begin{align}\nonumber
	\int_{\mathbb{R}^N}(\mathbf{I}_\alpha[\eta\otimes \delta_{\{x_{_N}=0\}}] )^{s'}x_ndx&\nonumber\geq 	\int_{\mathbb{R}^N}\left(\int_{2x_{_N}}^{4x_{_N}}\frac{(\eta\otimes \delta_{\{x_{_N}=0\}})(B_r(x))}{r^{N-\alpha}}\frac{dr}{r}\right)^{s'}x_ndx\\
	& \nonumber\geq c_{_{97}}	\int_{\mathbb{R}^N}\left(\frac{\eta(B'_{x_{_N}}(x'))}{x_{_N}^{N-\alpha}}\right)^{s'}x_{_N}dx\\ &\geq c_{_{98}} \int_{\mathbb{R}^{N-1}}\left(\sup_{t>0}\frac{\eta(B'_t(x'))}{t^{N-1-\alpha-\frac{2}{s'}}}\right)dx.'\label{es5}
	\end{align}
	By using Lemma \ref{hardy} we obtain
	\bel{es6}\BA{lll}\displaystyle
	\int_{\mathbb{R}^N}(\mathbf{I}_\alpha[\eta\otimes \delta_{\{x_{_N}=0\}}] )^{s'}x_ndx\leq  	\int_{\mathbb{R}^N}\left(\int_{x_{_N}}^{\infty}\frac{\eta(B'_r(x'))}{r^{N-\alpha}}\frac{dr}{r}\right)^{s'}dx_{_N}dx'\\[4mm]\phantom{\int_{\mathbb{R}^N}(\mathbf{I}_\alpha[\eta\otimes \delta_{\{x_{_N}=0\}}] )^{s'}x_{_N}dx}
	\displaystyle
	\leq c_{_{99}}\int_{\mathbb{R}^{N-1}}\int_{0}^{\infty}\left(\frac{\eta(B'_t(x'))}{t^{N-1-\alpha-\frac{2}{s'}}}\right)^{s'}\frac{dt}{t}dx'.
	\EA\ee
	On the other hand, by \cite[Proposition 5.1]{22PhVe}, there holds
\bel{es7}\BA{lll}\displaystyle
	\int_{\mathbb{R}^{N-1}}\left(\sup_{t>0}\frac{\eta(B'_t(x'))}{t^{N-1-\alpha-\frac{2}{s'}}}\right)dx'\asymp \int_{\mathbb{R}^{N-1}}\int_{0}^{\infty}\left(\frac{\eta(B'_t(x'))}{t^{N-1-\alpha-\frac{2}{s'}}}\right)^{s'}\frac{dt}{t}dx'\\[4mm]
	\phantom{\int_{\mathbb{R}^{N-1}}\left(\sup_{t>0}\frac{\eta(B'_t(x'))}{t^{n-1-(\alpha+\frac{2}{s'})}}\right)dx'}\displaystyle
	\asymp \int_{\mathbb{R}^{N-1}}\left(\int_{0}^{\infty}\frac{\eta(B'_t(x'))}{t^{N-1-\alpha-\frac{2}{s'}}}\frac{dt}{t}\right)^{s'}dx'.	
	\EA\ee
Combining \eqref{es5}, \eqref{es6} and \eqref{es7} we obtain \eqref{es8}. Moreover, we deduce \eqref{es9}  from \eqref{es8} and \cite[Theorem 2.5.1]{22AH}, which ends the proof.
\qeda
\medskip

\nind{\bf Proof of Theorem \ref{traceproblem}} The following estimates are cclassical 
\begin{align}
	& \mathbf{G}(x,y)\asymp \frac{x_{_N}y_{_N}}{|x-y|^{N-2}\max\{|x-y|,x_{_N},y_{_N}\}^2}\leq c_{_{100}} \frac{y_{_N}}{|x-y|^{N-1}},         \\&
	\mathbf{P}(x,z)= c_{_{101}}\frac{x_{_N}}{|x-z|^{N}}\leq c_{_{101}} \frac{1}{|x-z|^{N-1}}.
	\end{align}   
	Thus, 
	\begin{align}
	\mathbf{G}\left[(\mathbf{P}[\sigma_1])^{q_2}\right]+\mathbf{P}[\sigma_2]\leq c_{_{102}} \mathbf{I}_1[\omega],
	\end{align}   
	where $\omega(x)=\gr (\mathbf{P}[\sigma_1])^{q_2}+\sigma_2 $ in $\mathbb{R}^N$.  
	Therefore, we infer that if 
	\begin{align}\label{2412201412}
	\mathbf{I}_1\left[\left(\mathbf{I}_2\left[\left(\mathbf{I}_1[\omega]\right)^{q_1}\right]\right)^{q_2}\chi_{_{\mathbb{R}^N_+}}\rho\right]\leq c_{_{103}} \mathbf{I}_1[\omega]~~\text{ in }~\mathbb{R}^N_+
	\end{align}
	for some $c_{_{103}}>0$ small enough, then \eqref{0411201410d*} admits a positive solution $(u,v)$.  On the other hand, we deduce \eqref{2412201412} from  Lemma \ref{lela}  and Theorem \ref{final}. The proof is complete.
\qeda


\begin{remark} The system 
\begin{equation}\label{BVY}
\begin{array}
{llll}
-\Delta u=v^{q_1}+\ge_1\gm\qquad&\text{in} ~\Gw\\[1mm]
-\Delta v=u^{q_2}+\ge_2\eta\qquad&\text{in} ~\Gw\\[1mm]  
\phantom{-\Delta}
u=\ge_3\sigma_1,v=\ge_4\sigma_2\qquad&\text{in} ~\prt\Gw,                                                                               
\end{array}
\end{equation}
where $d(.)\gm,d(.)\gl$ belong to $\mathfrak M_b^+(\Gw)$, $\sigma_1,\sigma_2$ to $\mathfrak M^+(\prt\Gw)$ and the $\ge_j$ are positive numbers,
is analyzed in \cite[Th 4.6]{22Bixy}. Therein it is proved that if 
\begin{equation}\label{BVY-1}
\BA {lll}
\myint{\Gw}{}\left(\mathbf{G}[\gm]+\mathbf{P}[\gl]\right)^{\max\{q_1,q_2\}}d(x) dx<\infty,                                                            
\EA
\end{equation}
which is equivalent to a capacitary estimate, and
\begin{equation}\label{BVY-2}
\BA {lll}
\min\left\{q_2\myfrac{q_1+1}{q_2+1},q_1\myfrac{q_2+1}{q_1+1}\right\} <\myfrac{N+1}{N-1},                                                           
\EA
\end{equation}
and if the $\ge_j$ are small enough, then (\ref{BVY}) admits a positive solution. Now condition (\ref{BVY-2}) is a subcriticality assumption (for at least one of the two exponents $q_j$) since there is no condition on the boundary measures. 
\end{remark}
 

\begin{thebibliography}{99} 
 	
 	\bibitem{22AH} D. R. Adams, L.I. Heberg, {\em Function Spaces and Potential Theory},  Grundlehren der Mathematischen Wisenschaften {\bf 31}, Springer-Verlag (1999).
 	
 	\bibitem{22AdPi} D.R. Adams, M. Pierre, {\em  Capacitary strong type estimates in semilinear problems}, Ann. Inst. Fourier (Grenoble) {\bf 41} (1991), 117-135.
 	
 	\bibitem{22BaPi} P. Baras, M. Pierre, {\em Crit\`ere d'existence des solutions positives pour des \'equations semi-lin\'eaires non monotones}, Ann. Inst. H. Poincar\'e, Anal. Non Lin. {\bf  3} (1985), 185-212.
 	
 	\bibitem{22Bi1}  M. F. Bidaut-V\'eron, {\em Local and global behavior of solutions of quasilinear equations of Emden-Fowler type}, Arch. Ration. Mech. Anal. {\bf 107} (1989), 293-324.
 	
 	\bibitem{22Bi2}  M. F. Bidaut-V\'eron, {\em  Necessary conditions of existence for an elliptic equation with source term and measure data involving p-Laplacian}, in: Proc. 2001 Luminy Conf. on Quasilinear Elliptic and Parabolic Equations and Systems, Electron. J. Differ. Equ. Conf. {\bf  8} (2002), 23-34.
 	
 	\bibitem{22Bi3} M. F. Bidaut-V\'eron, {\em Removable singularities and existence for a quasilinear equation with absorption or source term and measure data}, Adv. Nonlinear Stud. {\bf  3} (2003), 25-63.
 	
 	\bibitem{22VHV} M. F. Bidaut-V\'eron, Q.-H. Nguyen, L. V\'eron, {\em Quasilinear Lane-Emden equations with absorption and measure data}, J. Math. Pures Appl. {\bf  102} (2014), 315Ð337. 
	
\bibitem{VHV2} M.F. Bidaut-V\'eron, Q.-H. Nguyen, L. V\'eron; {\em Quasilinear elliptic  equations  with source mixed term and  measure data,} preprint.  	

 	\bibitem{22Bi4}  M. F. Bidaut-V\'eron, S. Pohozaev,  {\em Nonexistence results and estimates for some nonlinear elliptic problems}, J. Anal. Math. {\bf  84} (2001), 1-49.
 	
 	\bibitem{22Bixy}  M. F. Bidaut-V\'eron, C. Yarur, {\em Semilinear elliptic equations and systems with measure data: existence and a priori estimates}, Adv. Diff. Equ. {\bf  7} (2002), 257-296.

 	\bibitem{22BiDe}  I. Birindelli, F. Demengel, {\em Some Liouville theorems for the p-Laplacian}, in: Proc. 2001 Luminy Conf. on Quasilinear Elliptic and Parabolic Equations and Systems, Electron. J. Differ. Equ. Conf. {\bf  8} (2002), 35-46.
 	
 	\bibitem{22DMOP} G. Dal Maso, F. Murat, L. Orsina, A. Prignet, {\em Renormalized solutions of elliptic equations with general measure data}, Ann. Sc. Norm. Sup. Pisa {\bf  28} (1999), 741-808.
 	
 	
 	
 	
 	\bibitem{22HeKiMa} J. Heinonen, T. Kilpelainen, O. Martio, {\em Nonlinear potential theory of degenerate elliptic equations}.  Unabridged republication of the 1993 original. Dover Publications, Inc., Mineola, NY, 2006. xii+404 pp.
 	
 	\bibitem{22KiMa1} T. Kilpelainen, J. Mal\'y, {\em Degenerate elliptic equation with measure data and nonlinear potentials}, Ann. Sc. Norm. Super. Pisa, Cl. Sci. {\bf  19} (1992), 591-613. 
 	
 	\bibitem{22KiMa2} T. Kilpelainen, J. Mal\'y, {\em The Wiener test and potential estimates for quasilinear elliptic equations}, Acta Math. {\bf  172} (1994), 137-161.
 	
 	\bibitem{22HoJa} P. Honzik, B. Jaye, {\em On the good-$\lambda$ inequality for nonlinear potentials}, Proc. Amer. Math. Soc. {\bf  140} (2012), 4167-4180.%
 	
 	\bibitem{22La}D. Labutin, {\em Potential estimates for a class of fully nonlinear elliptic equations}, Duke Math. J. {\bf 111}, 1-49,  (2002).
 	\bibitem{MV5} M. Marcus, L. V\'eron, {\em Nonlinear second order elliptic equations involving measures.} De Gruyter Series in Nonlinear Analysis and Applications {\bf  21}. De Gruyter, Berlin (2014).
	
 	\bibitem{22ReRa} M. M. Rao, Z. D. Ren, {\em Theory of Orlicz Spaces}, Textbooks in Pure and Applied Mathematics (1991).
 	
 	\bibitem {22PhVe}N. C. Phuc, I. E. Verbitsky, {\em Quasilinear and Hessian
 		equations of Lane-Emden type}, Ann. Math. {\bf  168}, 859-914 (2008).
 	
 	\bibitem{22PhVe2}N. C. Phuc, I. E. Verbitsky, {\em Singular quasilinear and Hessian equation and inequalities}, J. Funct. Anal. {\bf  256} (2009), 1875-1906. 
 	
 	
 	\bibitem{22SeZo} J. Serrin, H. Zou, {\em Cauchy-Liouville and universal boundedness theorems for quasilinear elliptic equations and inequalities}, Acta Math. {\bf  189} (2002) 79-142.
 	
 	%
 	\bibitem{22TW1} N. S. Trudinger and X.J.  Wang, {\em Hessian measures}, Topological Methods in Nonlinear Analysis {\bf  10} (1997), 225-239.
 	
 	\bibitem{22TW2} N. S. Trudinger, X.J. Wang, {\em Hessian measures II}, Annals of Mathematics {\bf  150}  (1999), 579-604.
 	
 	\bibitem{22TW3} N. S. Trudinger, X. J.  Wang, {\em Hessian measures III}, Journal of Functional Analysis {\bf  193} (2002), 1-23.
 	
 	\bibitem{22TW4}  N. S. Trudinger, X. J.  Wang, {\em On the weak continuity of elliptic operators and applications to potential theory}, Amer. J. Math. {\bf  124} (2002), 369-410.
 	
 	\bibitem{22V}L. V\'eron, {\em Elliptic equations involving measures}, in  Stationary Partial  Differential Equations, Handbook of Equations {\bf vol. I}, Elsevier B.V., pp. 593-712 (2004). %
	
	\bibitem{VeQ}L. V\'eron, {\em Local and global aspects of quasilinear degenerate elliptic equations. Quasilinear elliptic singular problems,} World Scientific Publishing Co. Pte. Ltd., Hackensack, NJ, xv+457 (2017).
	
 \end{thebibliography}
\end{document}